\documentclass{article}
\usepackage{arxiv}
\usepackage{float}
\usepackage{algorithm,algorithmic}
\usepackage{amsmath,amssymb}
\usepackage{amsthm,latexsym,float}
\usepackage[labelformat=simple]{subfig}

\usepackage{bm,multirow}        
\usepackage{graphicx}
\usepackage{xcolor}
\usepackage{placeins}
\usepackage{doi}
\graphicspath{{./fig/}}

\usepackage{hyperref}

\title{Reduced order modelling of nonlinear cross-diffusion systems}

\author{ \href{https://orcid.org/0000-0003-1037-5431}{\includegraphics[scale=0.06]{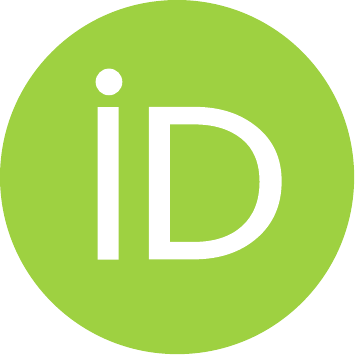}\hspace{1mm}B\"ulent Karas\"ozen} \\
     Institute of Applied Mathematics \& Department of Mathematics,
     Middle East Technical University\\
     Ankara-Turkey\\
     \texttt{bulent@metu.edu.tr}
	\And
\href{https://orcid.org/0000-0001-8952-7658}{\includegraphics[scale=0.06]{orcid.pdf}\hspace{1mm}G\"ulden M\"ulayim} \\
Department of Mathematics, Ad{\i}yaman University, Ad{\i}yaman-Turkey\\
\& Institute of Applied Mathematics, Middle East Technical University,  Ankara-Turkey\\
 \texttt{gulden.mulayim@metu.edu.tr}
	\And
  \href{https://orcid.org/0000-0001-5262-063X}{\includegraphics[scale=0.06]{orcid.pdf}\hspace{1mm}Murat Uzunca} \\
   Department of Mathematics, Sinop University,
     Sinop-Turkey \\
     \texttt{muzunca@sinop.edu.tr}\\
     \And
 \href{https://orcid.org/0000-0001-7904-605X}{\includegraphics[scale=0.06]{orcid.pdf}\hspace{1mm}S\"uleyman Y{\i}ld{\i}z} \\
Institute of Applied Mathematics\\
Middle East Technical University, 
 Ankara-Turkey\\
	\texttt{yildiz.suleyman@metu.edu.tr} \\
}

\date{}

\begin{document}
\maketitle

\begin{abstract}
In this work, we present a reduced-order model for a nonlinear cross-diffusion problem from population dynamics, for the Shigesada-Kawasaki-Teramoto (SKT) equation with Lotka-Volterra kinetics. The finite-difference discretization of the SKT equation in space results in a system of linear--quadratic ordinary differential equations (ODEs). The reduced order model (ROM) has the same linear-quadratic structure as the full order model (FOM). Using the linear-quadratic structure of the ROM, the reduced-order solutions  are computed independent of the full solutions with the proper orthogonal decomposition (POD). The computation of the reduced solutions is further accelerated by applying tensorial POD. The formation of the patterns of the SKT equation consists of a  fast transient phase and a long steady-state phase. Reduced  order solutions are computed by separating the time, into two-time intervals. In numerical experiments, we show for one-and two-dimensional SKT equations with pattern formation, the reduced-order solutions obtained in the time-windowed form, i.e., principal decomposition framework (P-POD),   are more accurate than   the global POD solutions (G-POD) obtained in the whole time interval.  Furthermore, we show the decrease of the entropy numerically by the reduced solutions, which is important for the global existence of nonlinear cross-diffusion equations such as the SKT equation.
\end{abstract}


\section{Introduction}

The interaction between species has been widely studied with reaction-diffusion models. Cross-diffusion systems are quasilinear parabolic equations in which the gradient of one variable induces a flux of another variable. They arise in multi-component systems from physics, chemistry, and biology. The correlation between diffusion and cross-diffusion terms may cause an unstable steady-state, called Turing instability or diffusion-driven instability, which leads to the formation of patterns \cite{Gambino12,Gambino13}. In this paper we consider a well-known cross-diffusion system from population dynamics, the Shigesada-Kawasaki-Teramoto (SKT) equation with Lotka-Volterra kinetics \cite{Shigesada79}. Numerical methods for coupled PDEs like the SKT equation are computationally expensive and require a large amount of computer memory and computing time in real-time simulations. During the last decades, reduced-order models (ROMs) have emerged as a powerful approach to reduce the cost of evaluating large systems of PDEs by constructing a low-dimensional linear subspace (reduced space) that approximately represents the solution to the system of PDEs with a significantly reduced computational cost. Semi-discretization of the SKT equation in space by finite-differences results in a linear-quadratic system. Recently reduced-order methods are developed for linear-quadratic systems; among them are moment-matching \cite{Ahmad17}, ${\mathcal H}_2$-quasi-optimal model order reduction \cite{Benner18}, balanced truncation \cite{Benner17a},  data-driven methods for quadratic-bilinear systems \cite{Antoulas18}. The proper orthogonal decomposition (POD) has been widely used as a computationally efficient reduced-order modeling technique in large-scale numerical simulations of nonlinear PDEs; the Navier-Stokes equation, Burger's equation, reaction-diffusion systems. The solutions of the high fidelity full-order model (FOM), generated by space-time discretization of PDEs, are projected onto a reduced space using the proper POD, resulting in a dynamical system with much lower order than the FOM \cite{Berkooz93,Sirovich87}. Applying POD Galerkin projection, dominant modes of the PDEs are extracted from the snapshots of the FOM solutions. The computation of the full-order solutions, POD basis functions, reduced matrices, and matricized tensors are performed in the "offline" stage, whereas the reduced system is solved in the "online" stage. We refer to the books \cite{Benner17,Quarteroni16} for an overview of the available ROM techniques.

There are many papers where, turbulent flow patterns, wave patterns, etc are investigated by various ROM methods. To the best of our knowledge, the formation of Turing patterns of the FitzHugh-Nagumo equation is investigated by POD-DEIM method only in \cite{Karasozen16,Karasozen17}. The SKT equation has linear and quadratic nonlinear terms, both in the cross-diffusion and in the Lotka-Volterra parts. Consequently, the semi-discretization of the SKT equation by second-order finite difference results in a system of linear-quadratic ordinary differential equations (ODEs). For time discretization, we use second-order linearly implicit Kahan's method \cite{Kahan97,Celledoni13}, which is designed for ODEs with quadratic nonlinear terms, as the SKT equation. In contrast to the fully implicit schemes, such as the Crank-Nicolson scheme, Kahan's method requires only one step Newton iteration at each time step.

The reduced system obtained by Galerkin projection contains nonlinear terms depending on the dimension of FOM, i.e., the offline and online phases are not separated. Using hyper-reduction techniques like discrete empirical interpolation method \cite{Chaturantabut12}, in the offline stage, additional reduced space is constructed approximating the nonlinear terms in the FOM, which may cause inaccuracies in the ROM solutions. When nonlinear PDEs like the SKT equation have polynomial structure, projecting the FOM onto the reduced space yields low-dimensional matrix operators that preserve the polynomial structure of the FOMs. {We apply POD by exploiting matricization of tensors \cite{Benner18,Benner15,Kramer19} to the quadratic terms of the semi-discrete SKT equation, so that the offline and online phases are separated. This enables the construction of computationally efficient ROMs without using hyper-reduction techniques like discrete empirical interpolation method. Here we make use of the sparse matrix technique MULTIPROD \cite{Leva08mmm} to speed up the matricized tensor calculations further.

For smooth systems where the system energy can be characterized by using a few modes, the global-POD (G-POD) method in the whole time interval provides a very efficient way to generate reduced-order systems. However, its applicability to complex, nonlinear PDEs is often limited by the errors associated with the finite truncation in POD modes and by the fast/slow solutions in different regimes. An alternative (and complementary) approach is the partitioned-POD (P-POD) which is based on the principal interval decomposition (PID) \cite{Borggaard15,Ahmed19,Ahmed20} that optimizes the length of time windows over which the POD procedure is performed. The cross-diffusion systems with pattern formation as the SKT equation have a rapidly changing short transient phase and long stationary phase. These two phases provide a natural decomposition of the whole time domain into two sub-intervals in the PID framework. We show for one- and two-dimensional SKT equations, the patterns can be more efficiently and accurately computed with P-POD than with the G-POD. It is critical to determining whether the ROM can issue reliable predictions in regimes outside of the training data. We show for one- and two-dimensional problems, the patterns can be predicted by the ROMs with acceptable accuracy. Moreover, cross-diffusion systems have an entropy structure. We show that the dissipation of reduced entropy is well preserved for the SKT equation.

The organization of this paper is as follows. In Section~\ref{sktsec}, we briefly describe the SKT equation. The fully discrete model in space and time is derived in Section~\ref{fomsec}. The reduced order methods G-POD, P-POD, and tensor techniques are described in Section~\ref{romsec}. Numerical experiments for one- and two-dimensional SKT equations are presented in Section~\ref{numsec}. Finally, we provide brief conclusions and directions for future work.

\section{Shigesada-Kawasaki-Teramoto equation}
\label{sktsec}

The interaction between species has been widely studied with reaction-diffusion models. The most prominent example is the Lotka-Volterra competition diffusion system which has been extensively studied in population ecology. When the diffusion of one of the species depends not only on the density of these species but also on the density of the other species, then cross-diffusion occurs, which may give rise to the formation of patterns. The species with high densities diffuse faster than predicted by the usual linear diffusion towards lower density areas, which leads to the coexistence of two spatial segregated competing species, known as cross-diffusion induced instability. When cross-and self-diffusion is absent, for linear diffusion in a convex domain, the only stable equilibrium solutions are spatially homogenous. In reaction-diffusion with cross-diffusion, the destabilization of a constant steady-state is followed by the transition to a non-homogeneous steady-state, i.e., the formation of patterns which are stationary in time and periodic or inhomogeneous in space. The linear stability analysis shows that the cross-diffusion is the key mechanism for the formation of spatial patterns through Turing instability \cite{Turing52}. In this paper, we consider the strongly coupled reaction-diffusion system with nonlinear self- and cross-diffusion terms. One of the most popular models in population ecology with pattern formation is the SKT cross-diffusion system \cite{Shigesada79} with the Lotka-Volterra reaction terms
\begin{equation}\label{skt}
\begin{aligned}
\frac{\partial u_1}{\partial t} &=\Delta (c_1  +a_1 u_1 +b_1 u_2) u_1+ \Gamma (r_1  - \gamma_{11}
u_1 - \gamma_{12} u_2) u_1  \\
\frac{\partial u_2}{\partial t} &= \Delta (c_2  +a_2 u_2 +b_2 u_1)u_2 +\Gamma (r_2 - \gamma_{21} u_1 - \gamma_{22} u_2)u_2 \\
\end{aligned}
\end{equation}
in a convex bounded domain $\Omega \subset \mathbb{R}^d,\; (d=1,2)$ with a smooth boundary $\partial \Omega$ on the time interval
$t\in[0,T] \subset \mathbb{R}$ with $T > 0$.
In Eq  \eqref{skt}, $u_1(x,t)$ and $u_2(x,t)$ with $x\in \Omega \subset \mathbb{R}^d$ denote population densities of two competing species and $\Delta $ is the Laplace operator. The initial and boundary conditions are
$$
u_1(x,0) =u_1^0(x), \ u_2(x,0) =u_2^0(x) \quad \text{in }  \Omega, \quad \frac{\partial u}{\partial {\bm n}}  = \frac{\partial v}{\partial {\bm n}}
= 0 \quad \text{on } \partial \Omega \times (0,T),
$$
where ${\bm n}$ is the unit outward normal vector to the boundary $\partial \Omega$. The homogeneous Neumann (zero-flux) boundary conditions impose the weakest constraint on formation of self-organizing patterns \cite{Gambino12,Gambino13}.

The parameters $a_i$ and $c_i$ are self-diffusion and linear diffusion coefficients, respectively, while the parameters $b_i$ are the cross-diffusion coefficients. The parameters $a_i, b_i, c_i, \gamma_{ij}, \; (i,j=1,2)$ are assumed to be non-negative. The parameters $r_i$ denote the intrinsic growth rates, $\gamma_{ii}$ the intra-specific competition coefficients, and $\gamma_{ij}, (i\ne j)$ are inter-specific competition rates. The parameter $\Gamma$ represents the relative strength of reaction terms.

Pattern formation in the SKT system \eqref{skt} was investigated using linear and weakly nonlinear stability analysis in \cite{Gambino12,Gambino13}. Cross-diffusion destabilizes the uniform equilibrium leading to traveling fronts \cite{Gambino12} in the one-dimensional SKT equation \eqref{skt} and formation of patterns in the two-dimensional SKT equation \eqref{skt} \cite{Gambino13}. In both papers, it was shown, for parameter values $b_1 > b_1^c$, patterns start to emerge from an initial condition with a random periodic perturbation of the equilibrium $(u_1^*, u_2^*)$
\begin{equation}\label{equilib}
(u_1^* ,u_2^*) = \left (  \frac{r_1 \gamma_{22} - r_2 \gamma_{12}}{\gamma_{11} \gamma_{22} -  \gamma_{12}\gamma_{21}},
 \frac{r_2 \gamma_{11} - r_1 \gamma_{21}}{\gamma_{11} \gamma_{22} -  \gamma_{12}\gamma_{21}} \right ).
\end{equation}
The critical value of bifurcation parameter $b_1^c$ is calculated using the Turing instability analysis for one- and two-dimensional SKT equations \eqref{skt} in \cite{Gambino12,Gambino13}.

The entropy structure is crucial to understand various theoretical properties of
cross-diffusion systems, such as existence, regularity and long time asymptotic
weak solutions of the SKT equation \eqref{skt}.
The SKT equation \eqref{skt} can be written alternatively \cite{Jungel15a,Jungel17,Chen19}
\begin{equation}\label{skt2}
\frac{\partial  u}{\partial t} = \text{div} \left ( {\mathcal A}(u) \nabla u\right) + f(u)
\end{equation}
with the diffusion matrix
$$
{\mathcal A} (u) =
\begin{pmatrix}
c_1 + 2a_1u_1 + b_1u_2 & b_1u_1 \\
b_2u_2 & c_2 + b_2u_1 + 2a_2u_2
\end{pmatrix},
$$
where $u=(u_1,u_2)^T$ and $f(u)=(f_1(u),f_2(u))^T$ with the Lotka-Volterra reaction terms $f_i(u)=\Gamma (r_i  - \gamma_{i1}u_1 - \gamma_{i2} u_2) u_i$, $i=1,2$. A characteristic feature of the cross-diffusion is that
the diffusion matrix is generally neither symmetric nor positive definite which complicates the mathematical analysis. However, using a transformation of variables (called entropy variable), the transformed diffusion
matrix becomes positive definite and sometimes even symmetric. Hence the existence of global solutions can be established. The entropy ${\mathcal E}$ for the SKT equation \eqref{skt2} is given without the reaction term $f(u)$ as \cite{Jungel15a,Jungel17,Chen19}
\begin{equation}\label{ent}
{\mathcal E}(u) = \int_{\Omega} h(u)\; dx, \quad h(u) =
\pi_1u_1(\log u_1 -1)  + \pi_2u_2(\log u_2 -1),
\end{equation}
when two constants $\pi_1$ and $\pi_2$ exist satisfying $\pi_1 b_1 = \pi_2 b_2$. For the given parameters $b_1$ and $b_2$, one can always find $\pi_1$ and $\pi_2$ satisfying this constraint. The entropy decreases in time, i.e., $\frac{d{\mathcal E}}{dt} \le 0$.

\section{Full order model}
\label{fomsec}

The SKT system \eqref{skt} has been solved by various numerical methods: fully implicit finite volume method \cite{Andreianov11}, semi-implicit finite difference method \cite{Galiano01}, semi-implicit spectral method \cite{Gambino13}, explicit Euler, and finite difference method \cite{Gambino12}. In \cite{ Murakawa11,Murakawa17}, the SKT system \eqref{skt} is transformed into a semi-linear PDE through replacing the nonlinear self- and cross-diffusion terms by linear reaction-diffusion terms. The resulting semi-linear equations with Lotka-Volterra reaction terms and linear diffusion terms are solved by the explicit Euler method with the finite difference or finite volume discretization in space.

Here we discretize the SKT equation \eqref{skt} in space by finite differences, which leads to a system of ODEs of the form
\begin{equation} \label{sktfd}
\begin{aligned}
	\frac{d \bm{u}_1}{d t} &= A(c_1\bm{u}_1 + a_1\bm{u}_1^2 + b_1 \bm{u}_1\odot\bm{u}_2 )  + \Gamma (r_1 \bm{u}_1 -
\gamma_{11} \bm{u}_1^2 - \gamma_{12}\bm{u}_1\odot \bm{u}_2),  \\
	\frac{d \bm{u}_2}{d t} &= A( c_2 \bm{u}_2 + a_2 \bm{u}_2^2  +b_2  \bm{u}_1\odot\bm{u}_2) + \Gamma (r_2\bm{u}_2 -
\gamma_{22}\bm{u}_2^2
-\gamma_{21}\bm{u}_1 \odot\bm{u}_2 ),
\end{aligned}
\end{equation}
where $\bm{u}_1(t),\bm{u}_2(t):[0,T]\to\mathbb{R}^N$ are semi-discrete approximations to the exact solutions $u_1(x,t)$ and $u_2(x,t)$ at $N$ spatial grid nodes, $\odot$ denotes the element-wise multiplication (Hadamard operator). The number $N$ of spatial grid nodes differ for one- and two-dimensional regions. The components of the semi-discrete solution vectors $\bm{u}_i(t)$ ($i=1,2$) in the case of one- and two-dimensional regions are given respectively by
\begin{align*}
\bm{u}_i(t) &= (u_{i}(x_1,t),\ldots,u_{i}(x_{n_x},t))^T, \\
\bm{u}_i(t) &= (u_{i}(x_1,y_1,t),\ldots,u_{i}(x_{n_x},y_1,t),u_{i}(x_1,y_2,t),\ldots,  (x_{n_x},y_{n_y},t))^T,
\end{align*}
with $N=n_x$ on one-dimensional regions and $N=n_xn_y$ on two-dimensional regions, where $n_x$ and $n_y$ are the number of partition in $x$ and $y$-directions, respectively.

In the ODE system \eqref{sktfd}, the matrix $A$ represents the finite difference matrix related to the second order centered finite differences approximation to the Laplace operator $\Delta$ under homogeneous Neumann boundary conditions. More clearly, let $I_n\in \mathbb{R}^{n\times n}$ denotes the $n$-dimensional identity matrix, and let the matrix $B_n\in \mathbb{R}^{n\times n}$ given by
\begin{equation}
B_n =
\begin{pmatrix}
-2 & 2 & & &  \\
1 & -2 & 1 & &  \\
& \ddots & \ddots & \ddots & \\
 & & 1 & -2 & 1\\
  & & & 2 & -2
\end{pmatrix},
\end{equation}
corresponding to the centered finite differences approximation to the Laplace operator under homogeneous Neumann boundary conditions with $n+1$ grid nodes. Then, on the one-dimensional domains, it is given as
$$
A = \frac{1}{\Delta x ^2}B_{n_x}\in \mathbb{R}^{n_x\times n_x},
$$
whereas on the two-dimensional domains we have
$$
A =  \frac{1}{\Delta x ^2}(B_{n_x}\otimes I_{n_y}) + \frac{1}{\Delta y ^2}(I_{n_x}\otimes B_{n_y}) \in \mathbb{R}^{n_xn_y\times n_xn_y},
$$
where $\Delta x$ and $\Delta y$ are the mesh sizes in $x$- and $y$-directions, respectively, and  $\otimes$ denotes the Kronecker product.

Collecting linear and quadratic parts, the system \eqref{sktfd} can be written as the following linear-quadratic ODE system
\begin{equation}
\begin{aligned}
	\frac{d \bm{u}_1}{d t} &=&   \underbrace{L_1 \bm{u}_1}_\text{linear}  +  \underbrace{Q_{11} \bm{u}_1^2 + Q_{12}(\bm{u}_1 \odot\bm{u}_2)}_\text{quadratic}, \\
	\frac{d \bm{u}_2}{d t} &=& \underbrace{L_2 \bm{u}_2}_\text{linear}  +
\underbrace{Q_{22}\bm{u}_2^2  + Q_{21}(\bm{u}_1\odot\bm{u}_2) }_\text{quadratic},
\end{aligned}
\end{equation}
where we set
$$
L_1 = c_1A + \Gamma r_1, \quad Q_{11} = a_1A - \Gamma \gamma_{11}, \quad Q_{12} = b_1A - \Gamma \gamma_{12},
$$
$$
L_2 = c_2A + \Gamma r_2, \quad Q_{21} = b_2A - \Gamma \gamma_{21}, \quad Q_{22} = a_2A - \Gamma \gamma_{22}.
$$
In compact form, we can also write as
\begin{equation}\label{linquad}
\frac{d \bm{u} }{dt } = \bm{F}(\bm{u}) =  L\bm{u}  +  M(\bm{u} \otimes \bm{u}) ,
\end{equation}
 where $\bm{u}=(\bm{u}_1,\bm{u}_2)   \in \mathbb{R}^{2N}$ is the state vector, $ L  \in \mathbb{R}^{2N\times 2N}$ represents the matrix of linear terms, and $M(\bm{u} \otimes \bm{u})= R_1(\bm{u}) + R_2(\bm{u}) $ is the quadratic part with  $R_i  \in \mathbb{R}^{2N}$ ($i=1,2$), given by
\begin{equation}\label{linfom}
\begin{aligned}
L&=
\begin{pmatrix}
L_1 & 0 \\
0 & L_2
\end{pmatrix},\\  R_1 &=
\begin{pmatrix}
Q_{11} & 0 \\
0 & Q_{22}
\end{pmatrix}
\begin{pmatrix}
H(\bm{u}_1\otimes\bm{u}_1)\\
H(\bm{u}_2\otimes\bm{u}_2)
\end{pmatrix}\;
R_2 =
\begin{pmatrix}
0 & Q_{12} \\
Q_{21} & 0
\end{pmatrix}
\begin{pmatrix}
H(\bm{u}_1\otimes\bm{u}_2)\\
H(\bm{u}_2\otimes\bm{u}_1)
\end{pmatrix}.
\end{aligned}
\end{equation}
In \eqref{linfom}, $H \in \mathbb{R}^{N\times N^2}$ stands for the matricized tensor so that it satisfies the identity $H(\bm{w}\otimes\bm{v})=\bm{w}\odot\bm{v}$ for any vector $\bm{w},\bm{v}\in \mathbb{R}^{N}$.

We solve the semi-discrete linear-quadratic ODE system \eqref{linquad} in time by Kahan's method \cite{Kahan93,Kahan97}:
\begin{equation}\label{kahan}
\frac{{\bm u}^{n+1} - {\bm u}^n}{\Delta t} =   \frac{1}{2}L({\bm u}^n + { \bm u}^{n+1}) + \tilde{R}_1({ \bm u}^n,{\bm u}^{n+1})+ \tilde{R}_2({ \bm u}^n,{\bm u}^{n+1}),
\end{equation}
where $\tilde{R}_i({ \bm u}^n,{\bm u}^{n+1})$ are the symmetric bilinear forms obtained by the polarization of  the quadratic vector fields $R_i$ \cite{Celledoni15}:
\begin{equation}
\tilde{R}_i({\bm  u}^n,{ \bm u}^{n+1}) = \frac{1}{2}R_i({ \bm  u}^n+{\bm u}^{n+1}) - R_i({\bm u}^n) -  R_i({\bm u}^{n+1}), \qquad i=1,2,
\end{equation}
$\Delta t$ is the time step size, and ${\bm u}^{n}$  is the full discrete solution vector at time $t_n$.
Kahan's method is time-reversal, symmetric and linearly implicit \cite{Celledoni13} for linear-quadratic systems such as the semi-discrete SKT equation
\eqref{linquad},  i.e., ${\bm u}^{n+1}$ can be computed by solving a single linear system
\begin{equation}
\left ( I_{2N} -\frac{\Delta t}{2} \bm{F}_J({\bm u}^n)\right ) \tilde{ {\bm u}} = \Delta t
\bm{F}({\bm u}^n),\qquad {\bm u}^{n+1} = {\bm u}^n +
\tilde{{\bm u}},
\end{equation}
where $\bm{F}_J({\bm u}^n)$ is the Jacobian matrix of $\bm{F}({\bm u})$ evaluated at ${\bm u}^n$. Kahan's method can also be written as a second order convergent Runge-Kutta method of the form \cite{Celledoni13}
\begin{equation} 
\frac{\bm{u}^{n+1} - \bm{u}^n}{\Delta t} = -\frac{1}{2}\bm{F}(\bm{u}^n) + 2\bm{F} \left (\frac{\bm{u}^{n+1} + \bm{u}^n}{2}  \right )   - \frac{1}{2}\bm{F}(\bm{u}^{n+1}).
\end{equation}
The linearly implicit Kahan's method is much faster than the fully implicit time integrators like the implicit Euler and mid-point rule. By linearly implicit we mean integrators that require the solution of precisely one linear system of equations in each time step. This is opposed to fully implicit schemes in which typically an iterative solver is applied that it may require the solution of several linear systems in each time step.

\section{Reduced order model}
\label{romsec}

In this section, we introduce ROMs for the SKT equation \eqref{skt}. The standard way of constructing ROMs is the use of POD with Galerkin projection on the whole time interval, G-POD. The solutions of the SKT equation \eqref{skt} converge quickly to steady-state after a short transient phase. We have constructed P-POD in two sub-intervals  respecting different behaviors of the average densities of the FOM solutions in the  transient and steady-state phases. In both approaches, the computation of the reduced quadratic nonlinear terms in the reduced system scales with the dimension of FOM. The quadratic terms  of the semi-discrete SKT equation \eqref{linquad} are efficiently evaluated using tensor techniques  so that the computational cost  of the reduced order solutions depends only on the dimension of the ROMs.

\subsection{Global-POD (G-POD)}

The POD basis vectors are computed using the method of snapshots. We refer to \cite{Hinze19,Kunisch01} for a detailed description of the POD with Galerkin projection. The POD basis for the semi-discrete SKT system \eqref{sktfd} can be obtained  by stacking all $\bm{u}_1$ and $\bm{u}_2$ in one vector $\bm{u}$ and to determine the  common subspace $V$ by taking the singular value decomposition (SVD) of that data. But this may produce unstable ROMs such that the resulting ROMs do not preserve the coupling structure of the PDE \cite{Benner15,Reiss07}. In order to maintain the coupling structure in ROMs of the coupled SKT equation, we compute the snapshot matrices and the POD basis vectors separately for the state components $\bm{u}_1$  and $\bm{u}_2$.
Consider the discrete state vectors $\bm{u}_1$ and $\bm{u}_2$ of  \eqref{sktfd}. The snapshot matrix $S_i$ corresponding to the state $\bm{u}_i$ is defined as
$$
S_i := [ \bm{u}_i^1, \cdots , \bm{u}_i^{N_t} ] \in\mathbb{R}^{N\times N_t}, \qquad i=1,2,
$$
where each column vector $\bm{u}_i^n\in\mathbb{R}^{N}$ is the full discrete solution vector at the time instance $t_n$, $n=1,\ldots ,N_t$. We then expand the SVD of the snapshot matrix $S_i$ with the singular values $\sigma_{i,1} \ge \sigma_{i,2} \ge \cdots\ge \sigma_{i,N_t}\ge 0$.

For a positive integer $k_i\ll N$, the $k_i$-POD basis matrix $V_{i,k_i} \in \mathbb{R}^{N\times k_i}$ minimizes the least squares error of the snapshot reconstruction by solving the following optimization problem
\begin{equation} \label{podopt}
\min_{V_{i,k_i} \in \mathbb{R}^{N\times k_i}} \|S_i  -V_{i,k_i}V_{i,k_i}^T S_i \|^2_F=
\min_{V_{i,k_i} \in \mathbb{R}^{N\times k_i}}\sum_{n=1}^{N_t} ||\bm{u}_i^n  -V_{i,k_i}V_{i,k_i}^T \bm{u}^n ||^2_2 = \sum_{n=k_i+1}^{N_t} \sigma_{i,n}^2,
\end{equation}
where $\|\cdot \|_2$ denotes the Euclidean $2$-norm, and $\|\cdot \|_F$ denotes the Frobenius norm. The optimal solution of basis matrix $V_{i,k_i}$ to this problem is given by the $k_i$ left singular vectors of $S_i$  corresponding to the $k_i$ largest singular values. The number of POD modes $k_i$ for each component $\bm{u}_i$ ($i=1,2$) is determined  usually by the relative information content  (RIC) defined by
\begin{equation} \label{ric}
\frac{\sum_{n=1}^{k_i} \sigma_{i,n}^2}{\sum_{n=1}^{N_t} \sigma_{i,n}^2  } < \text{tol}_\text{RIC}
\end{equation}
with a prescribed tolerance $\text{tol}_\text{RIC}$.

The POD state approximation is then given by $\bm{u}_i \approx V_{i,k_i} \widehat{\bm{u}}_i$, where $\widehat{\bm{u}}_i\in \mathbb{R}^{k_i}$ is the reduced state vector. Throughout the paper, we omit the subscript $k_i$ to simplify the notation and we denote just $V_i$ the $k_i$-POD basis matrix corresponding to the state $\bm{u}_i$.

Once the POD basis matrices $V_i\in\mathbb{R}^{N\times k_i}$ ($i=1,2$) are found, the ROM for the SKT system is obtained as the following linear-quadratic ODE system
\begin{equation} \label{romlinquad}
\frac{d }{dt }\widehat{\bm{u}} = \widehat{L}\widehat{\bm{u}}  + \widehat{R}_1 (\widehat{\bm{u}} )+ \widehat{R}_2 (\widehat{\bm{u}} ),
\end{equation}
where $\widehat{\bm{u}}=(\widehat{\bm{u}}_1,\widehat{\bm{u}}_2)$, $\widehat{L}=V^T LV$, $\widehat{R}_i(\widehat{\bm{u}}) = V^TR_i(V\widehat{\bm{u}} )$, and
$$
V = \begin{pmatrix}
V_1 & 0 \\
0 & V_2
\end{pmatrix} \in\mathbb{R}^{2N\times (k_1+k_2)}.
$$
The linear-quadratic structure of the FOM are preserved by the reduced model obtained by Galerkin projection \cite{Benner15a}.

\subsection{Partitioned POD (P-POD) }

POD depends on a global approximation of the snapshot data, which can result in overall deformation of the obtained modes for systems with fast variations in state and the constructed POD cannot capture any dominant structure at all. The P-POD approach is developed \cite{IJzerman01,Borggaard07,Borggaard15,Ahmed18,Ahmed19,Ahmed20} as an alternative approach that preserves the optimality of POD respecting local characteristics of the solutions. In P-POD, the time domain is divided into non-overlapping intervals, each characterizing a specific stage in the dynamics and evolution of the system. The same POD algorithm is applied within each sub-interval to generate a set of basis functions that best fit the respective partition locally. In short, the P-POD approach can be viewed as a decomposition of the G-POD subspace into a few locally optimal subspaces to obtain accurate partitioned ROMs with smaller sizes in each individual sub-interval. The P-POD was first studied in \cite{IJzerman01} and then is applied successfully to nonlinear convective fluid problems: to Burger's equation, Boussinesq equation, and Navier-Stokes equation \cite{Borggaard07,Borggaard15,Ahmed18,Ahmed19}, and to two-dimensional turbulence flow \cite{Ahmed20,Wang12}. Usually the time domain is decomposed into equidistant sub-intervals \cite{Ahmed18,Ahmed19,Ahmed20,Borggaard07}. Adaptive partitioning and clustering techniques can be used to construct the sub-intervals of different sizes in the time domain \cite{Borggaard15,Wang12}.

The cross-diffusion with pattern formation like SKT equation \eqref{skt} is characterized by a short transient phase and long stationary phase. In the transition phase, the solutions are changing rapidly. During the stationary phase, the solutions are changing slowly, until they reach the spatially inhomogeneous Turing patterns. This leads to the natural decomposition of the whole time domain into two sub-intervals in the P-POD approach, one for the transient phase and the other for the stationary phase. The transition from the transient phase to the stationary phase can be determined by the use of average densities
\begin{equation} \label{density}
\langle u_1(x,y,t)\rangle = \frac{1}{|\Omega|} \int_\Omega u_1(x,y,t) dxdy,  \quad \langle u_2(x,y,t)\rangle = \frac{1}{|\Omega|}  \int_\Omega u_2(x,y,t) dxdy,
\end{equation}
When the difference of the both average densities at two consecutive time instances are lower than a prespecified tolerance $\text{tol}_{\text{PID}}$, the stationary phase occur. This transition point is then used as the interface of two sub-intervals in the P-POD approach. More clearly, let $t_p$ denotes the transition point, i.e., the index $1<p<N_t$ is the minimum integer such that
$$
|\langle u_i(x,y,t_p)\rangle - \langle u_i(x,y,t_{p-1})\rangle| <\text{tol}_{\text{PID}} , \qquad i=1,2,
$$
for a given PID tolerance $\text{tol}_{\text{PID}}$. Then, we decompose the whole time-interval $[t_0,t_{N_t}]$ into two sub-intervals $I_1 = [t_0,t_p]$ and $I_2 = [t_p,t_{N_t}]$ with the common interface $t_p$. According to the P-POD approach, we set different POD basis matrices $V_i^{(1)}$ and $V_i^{(2)}$ through the snapshot matrices composed of the full solution vectors on either intervals $I_1$ and $I_2$,  respectively. Finally, we should enforce the interface constraint that the full discrete solution vectors $\bm{u}_i^{p}$ from the first interval $I_1$ agree with the initial vectors $\bm{u}_i^{p}$ on the second interval $I_2$. Using the POD basis matrices defined on two sub-intervals, this requires the following identity
$$
V_{i}^{(1)} \widehat{\bm{u}}_i^{(1)}(t_p) = V_{i}^{(2)} \widehat{\bm{u}}_i^{(2)}(t_p), \qquad i=1,2,
$$
where $\widehat{\bm{u}}_i^{(1)}(t_p)$ and $\widehat{\bm{u}}_i^{(2)}(t_p)$ are reduced solution vectors at the interface time $t_p$ from the intervals $I_1$ and $I_2$, respectively. Using the orthonormality of the POD basis matrices, the initial reduced solution vector on the interval $I_2$ can be recovered from the reduced solution vector at the interface $t_p$ from the interval $I_1$ as
$$
 \widehat{\bm{u}}_i^{(2)}(t_p) = (V_{i}^{(2)})^TV_{i}^{(1)} \widehat{\bm{u}}_i^{(1)}(t_p), \qquad i=1,2 .
$$

\subsection{ Efficient evaluation the nonlinear terms of the ROMs}

Although the dimension of the ROM \eqref{romlinquad} is small compared to the dimension of the FOM \eqref{sktfd} ($k_1+k_2\ll 2N$), the computation of the quadratic nonlinearities still scales with the dimension of FOM. This is overcome by   utilizing the matricized tensor together with the properties of  the  Kronecker product, and the offline and online stages are separated.

The explicit form of the reduced quadratic terms in the SKT equation \eqref{romlinquad} are given by
\begin{equation}
\begin{aligned}
\widehat{R}_1(\widehat{\bm{u}}) &= V^TR_1(V\widehat{\bm{u}} ) =
\begin{pmatrix}
V_1^TQ_{11}H((V_1\widehat{\bm{u}}_1)\otimes(V_1\widehat{\bm{u}}_1) ) \\
V_2^TQ_{22}H((V_2\widehat{\bm{u}}_2)\otimes(V_2\widehat{\bm{u}}_2) )
\end{pmatrix},\\
\widehat{R}_2(\widehat{\bm{u}}) &= V^TR_2(V\widehat{\bm{u}} ) =
\begin{pmatrix}
V_1^TQ_{12}H((V_1\widehat{\bm{u}}_1)\otimes(V_2\widehat{\bm{u}}_2) ) \\
V_2^TQ_{21}H((V_2\widehat{\bm{u}}_2)\otimes(V_1\widehat{\bm{u}}_1) )
\end{pmatrix}.
\end{aligned}
\end{equation}
It is clear that all the terms above are of the form
\begin{equation}\label{redtensor}
V_i^TQ_{ij}H((V_i\widehat{\bm{u}}_i)\otimes(V_j\widehat{\bm{u}}_j)), \qquad i,j=1,2.
\end{equation}

The terms with matricized tensor $H \in \mathbb{R}^{N\times N^2}$ in \eqref{redtensor} are computed using the properties of Kronecker product, which depends on the computation of the reduced matricized tensor $\widehat{H}= H(V_i\otimes V_j)\in\mathbb{R}^{N\times k_ik_j}$ so that we get
\begin{equation}\label{seperate}
V_i^TQ_{ij}H((V_i\widehat{\bm{u}}_i)\otimes(V_j\widehat{\bm{u}}_j)) = V_i^TQ_{ij}\widehat{H}(\widehat{\bm{u}}_i\otimes\widehat{\bm{u}}_j),
\end{equation}
where the small matrix $V_i^TQ_{ij}\widehat{H}\in\mathbb{R}^{k_i\times k_ik_j}$ can be precomputed in the offline stage. By applying tensor techniques to $ H(\bm{w}\otimes\bm{v})=\bm{w}\odot\bm{v} $, the high dimensional variables in \eqref{seperate} are  separated. Then the computational cost of the reduced quadratic nonlinear terms in the online stage depends only on the dimension of ROMs, i.e.  $\mathcal{O}(k_i^3)$. When  the  tensor techniques are not used, the online computation depends on the high dimensional FOM, e.g., $ (V_i\widehat{\bm{u}}_i)\odot(V_j\widehat{\bm{u}}_j) $, then the computational cost scales with the  dimension $N$  of the FOM, i.e.   $\mathcal{O}(Nk_i)$ \cite{Stefanescu14}.

Although $\widehat{H}$ is computed in the offline stage, the explicit computation of $V_i \otimes V_j$ may be inefficient since it depends on the full dimension $N$. In order to avoid from this computational burden,    $V_i \otimes V_j$ is computed in an efficient way using $\mu$-mode matricization of tensors \cite{Benner15t}.
Recently tensorial algorithms are developed by exploiting the particular structure of Kronecker product \cite{Benner18,Benner21}. The reduced matrix $\widehat{H}$ can be given in MATLAB notation as follows	
	\begin{align}\label{goyal}
	\widehat{H} =
	\begin{pmatrix}
	V_i(1,:)\otimes V_j(1,:)\\
	\vdots\\
	V_i(N,:)\otimes V_j(N,:)
	\end{pmatrix},
	\end{align}
which utilizes the structure of $ H(V_i\otimes V_j) $ without constructing $H$ explicitly. In \cite{Benner18,Benner21} the CUR matrix approximation \cite{Mahoney09} of $H(V_i\otimes V_j)$ is used to increase computational efficiency. Instead, here we make use of the "MULTIPROD" \cite{Leva08mmm} to increase computational efficiency. The MULTIPROD \footnote{\href{https://www.mathworks.com/matlabcentral/fileexchange/8773-multiple-matrix-multiplications-with-array-expansion-enabled}{https://www.mathworks.com/matlabcentral/fileexchange/8773-multiple-matrix-multiplications-with-array-expansion-enabled}} handles multiple multiplications of the multi-dimensional arrays via virtual array expansion. It is a fast and memory efficient generalization for arrays of the MATLAB matrix multiplication operator. For any given two vectors $\mathbf{a}$ and $\mathbf{b}$, the Kronecker product satisfies
	\begin{equation}\label{vec}
	\begin{aligned}
	(\text{vec}{(\mathbf{b}\mathbf{a}^\top)})^\top =(\mathbf{a}\otimes \mathbf{b})^\top
	=\mathbf{a}^\top\otimes \mathbf{b}^\top,
	\end{aligned}
	\end{equation}
where vec$(\cdot)$ denotes the vectorization of a matrix. Using the identity in \eqref{vec}, the matrix $\widehat{H}= H(V_i\otimes V_j)$ can be constructed as
\begin{equation}
\widehat{H}(m,:)=(\text{vec}(V_i(m,:)^\top V_j(m,:))^\top , \ \ m\in\{1,2,\ldots,N\}.
\end{equation}
Reshaping the matrix $V_i\in \mathbb{R}^{N\times k_i}$ as $\widetilde{V}_i \in \mathbb{R}^{N\times 1 \times k_i}$ and computing MULTIPROD of $V_j$ and $\widetilde{V}_i$ in $2$nd and $3$th dimensions, we obtain
$$
\mathcal{\widehat{H}} =\text{MULTIPROD}( V_j,\widetilde{V}_i)\in \mathbb{R}^{N\times k_j \times k_i},
$$
where the reduced matricized tensor $\widehat{H}\in\mathbb{R}^{N\times k_ik_j}$ is taken by reshaping the matrix $\mathcal{\widehat{H}}$ into a matrix of dimension $N\times k_ik_j$. The computation of $\mathcal{\widehat{H}}$ in \eqref{goyal} requires $N$ for loops within each iteration the matrix product of two matrices of sizes $k_j\times 1$ and $1\times k_i$ are done. But with the MULTIPROD, the matrix products are computed simultaneously in a single loop, decreasing the computational cost for the calculation of the reduced matricized tensor in the offline stage \cite{karasozen21}.

Alternatively, matricized tensor $\widehat{H}\in\mathbb{R}^{N\times k_ik_j}$ can be computed using the definition $H(\bm{w}\otimes\bm{v})=\bm{w}\odot\bm{v} $ via Hadamard product over reduced dimensions $ k_i $ and $ k_j $ as in Algorithm \eqref{algo1}. We remark that  in  the Algorithm \eqref{algo1} matricized tensor  is not evaluated over the FOM dimension $N$,  instead over the reduced dimensions  $ k_i $ and $ k_j $, which are small in our case.  
\begin{algorithm}
	\caption{Construction of matricized tensor $\widehat{H}\in\mathbb{R}^{N\times k_ik_j}$ over the reduced dimensions.}
	\begin{algorithmic}\label{algo1}
	\STATE\textbf{Input:} POD basis matrices $ V_i $ and $ V_j $\\
	\STATE \textbf{Output:} Matricized tensor $ \widehat{H} $
	\STATE Set $ k=1 $
    \FOR{$ m\in\{1,\ldots,k_i\} $}
    \FOR{$ n\in\{1,\ldots,k_j\} $}
    \STATE $ \widehat{H}(:,k) = V_i(:,m) \odot V_j(:,n) $
    \STATE $ k=k+1 $
    \ENDFOR
    \ENDFOR
	\end{algorithmic}
\end{algorithm}

\FloatBarrier
\section{Numerical results}
\label{numsec}

In this section we present results of the  numerical experiments for the one- and two-dimensional SKT system \eqref{skt}. We  compare the FOM and ROM solutions by G-POD and P-POD, and show the decreasing structure of the entropy. Furthermore we show  prediction capabilities of the ROMs in regimes outside of the training data. Recently similar predictions have been employed for combustion problems using Galerkin ROM \cite{Huang19} and data-driven ROM techniques \cite{Swischuk20}.

For examples in Section~\ref{exmp1} and Section~\ref{exmp2}, the initial conditions are taken as random periodic perturbation around the equilibrium $(u_1^*, u_2^*)$ given in \eqref{equilib} using MATLAB function {\bf rand\/},  uniformly distributed pseudo-random numbers. The reduced systems are computed with the same time step and the time integrator as the FOM, i.e., Kahan's method.

We expect the system to reach a spatially inhomogeneous steady-state. Thus, the simulation is left to run until the following criteria is satisfied
\begin{equation}\label{crit}
\frac{||\bm{u}_i^n - \bm{u}_i^{n-1} ||_{L^2(\Omega)}}{||\bm{u}^n||_{L^2(\Omega)}} \le \text{tol}_\text{ST}, \qquad i=1,2,
\end{equation}
where $\|\cdot\|_{L^2(\Omega)}$ denotes the usual $L^2$-norm over the domain $\Omega$, and calculated by the trapezoidal quadrature. The quantity \eqref{crit} is related to the rate of change of the variables, and so we simply stop the calculations when the solutions stop varying significantly with time for a prescribed tolerance $\text{tol}_\text{ST}$ $> 0$. We take in all simulations $\text{tol}_\text{ST}=10^{-6}$.

The accuracy of the ROM solutions are measured using the time averaged relative $L^2$-error defined by
\begin{equation}\label{relerr}
\|\bm{u}-V\widehat{\bm{u}}\|_{\text{rel}}=\frac{1}{N_t}\sum_{n=1}^{N_t}\frac{\|\bm{u}^n-V\widehat{\bm{u}}^n\|_{L^2(\Omega)}}{\|\bm{u}^n\|_{L^2(\Omega)}}.
\end{equation}

All the simulations are performed on a machine with Intel CoreTM i7 2.5 GHz 64 bit CPU, 8 GB RAM, Windows 10, using 64 bit MatLab R2014. For the time-dependent problems with many time steps, such as the SKT system, the snapshot matrix is large, leading to an expensive SVD. We use randomized SVD (rSVD) algorithm \cite{Halko11a} which only needs to perform SVD of small matrices, to efficiently generate a reduced basis with large snapshot matrices. In the ROM computations, tensor techniques are applied in the G-POD and P-POD.

\subsection{One-dimensional SKT equation  }
\label{exmp1}

Our first example is the one-dimensional SKT equation \eqref{skt} with
the parameters are taken from \cite{Gambino12}
\begin{align*}
a_1 &= 0.0001,\; a_2 = 0.1, \; b_1 = 6.5, \; b_2 = 0.3, \; c_1 = 0.2,\; c_2 = 0.2, \\
 \Gamma &= 49.75, \; r_1 =  1.2, \; r_2=1, \;\gamma_{11} = 0.5, \; \gamma_{12} = 0.4,\; \gamma_{21} = 0.38,\; \gamma_{22} = 0.41,
\end{align*}
where,  $b_1$ is taken larger than the critical value of the bifurcation parameter $b_c = 5.297$, so that pattern formation can occur. Spatial interval is set to $\Omega=[-\pi,\pi]$ with the mesh size $\Delta x =2\pi/200$ and $n_x =200$. The time step size is taken as $\Delta t = 0.001$. The steady-states are reached at $t=11.219$.

In Figure~\ref{1Dden}, left, the patterns at the steady-state are shown, that are formed starting from an initial datum which is a random periodic perturbation of the equilibrium \eqref{equilib}. The FOM solutions are very close to those in \cite{Gambino12}. In Figure~\ref{1Dden}, right, the densities start a plateau around $t=5$. Accordingly, using the PID tolerance $\text{tol}_{\text{PID}}=10^{-8}$, we obtain the transition point $t_p=4.023$, and the time interval is split into sub-intervals $I_1=[0,4.023]$ and $I_2=[4.023,11.219]$.

\begin{figure}[H]
\centering
\includegraphics[width=0.9\linewidth]{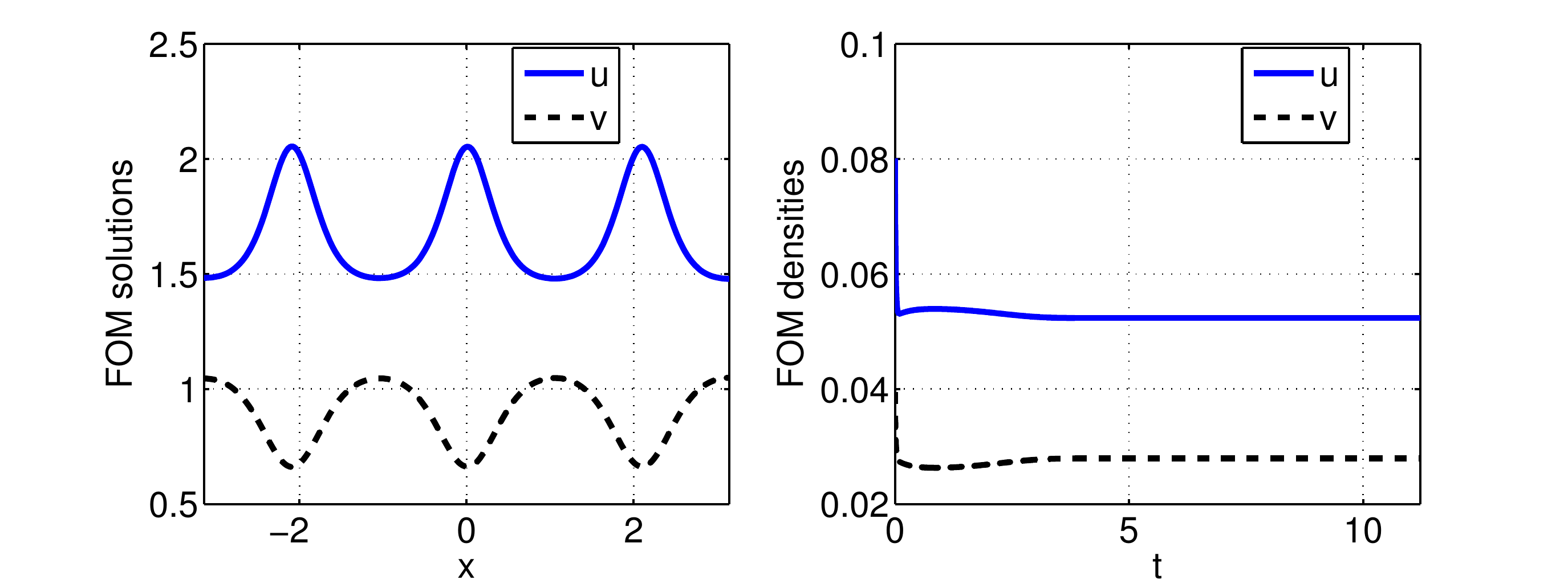}
\caption{FOM solutions \& densities for one-dimensional problem }
\label{1Dden}
\end{figure}

In Figure~\ref{1Dsvd}, normalized singular values $\sigma_{ij}/\sigma_{i1}$ are plotted on the whole time interval, on the intervals $I_1$ of the transient and $I_2$ of the steady-state phases, $i=1,2$. The singular values decay at the same rate, slowly on the whole time interval and on the first interval $I_1$, whereas the decay is faster on the second interval $I_2$ of the steady-state phase.

\begin{figure}[H]
\centering
\includegraphics[width=\linewidth]{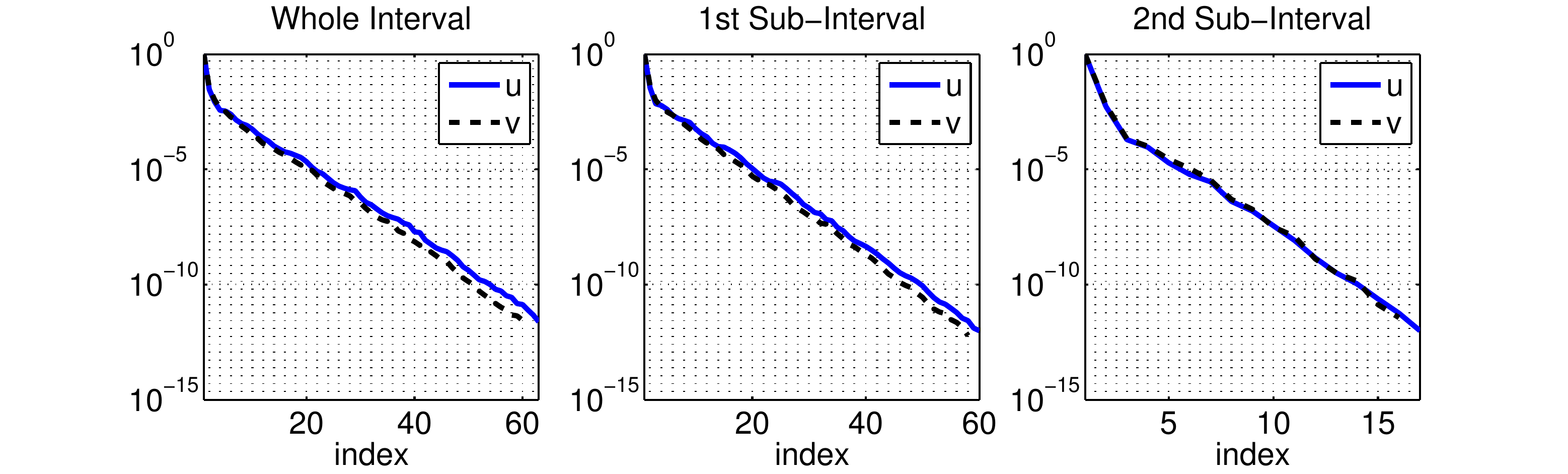}
\caption{Decay of normalized singular values for one-dimensional problem}
\label{1Dsvd}
\end{figure}

The number of POD modes and the time averaged relative $L^2$-errors \eqref{relerr} between FOM and ROM approximations for different RIC tolerances $\text{tol}_\text{RIC}$ in \eqref{ric} are listed in Table \ref{tab1D}. For the same RIC tolerances, the P-POD requires fewer POD modes on the interval $I_2$ of the steady-states comparing with the ones required on the interval $I_1$ of the transition phase. The time averaged relative $L^2$-errors obtained by the P-POD are smaller than the errors obtained by the G-POD approach. The ROM solutions in Figure~\ref{1Dsol} computed using the tolerance $\text{tol}_\text{RIC} = 10^{-4}$  are very close to the FOM solutions. To show the time evolution  of the patterns by FOM and ROMs, a video component is available and accompanies the electronic version of this paper. To access this video component, simply click on the image in Figure~\ref{1Dsol} (online version only).

\begin{table}[H]
  \centering
	\caption{Number of POD modes and errors for one-dimensional problem \label{tab1D}}
\resizebox{\columnwidth}{!}{
\begin{tabular}{|c|ll|lll|}\hline
\multirow{2}{*}{$\text{tol}_\text{RIC}$} & \multicolumn{2}{c|}{G-POD $u(v)$} & \multicolumn{3}{c|}{P-POD $u(v)$} \\
 &  \#Modes  & Error   & \#Modes-$I_1$  & \#Modes-$I_2$  & Error  \\
\hline
$10^{-3}$ & 5(5) & 5.35e-03(7.20e-03) & 6(5) & 2(2) & 1.84e-03(2.35e-03) \\
$10^{-4}$ & 9(7) & 9.23e-04(1.15e-03) & 9(8) & 2(2) & 7.74e-04(9.16e-04) \\
$10^{-5}$ & 12(10) & 4.14e-04(5.01e-04) & 12(10) & 2(2) & 2.15e-04(2.62e-04) \\
$10^{-6}$ & 15(14) & 1.34e-04(1.64e-04) & 15(13) & 3(4) & 9.27e-05(1.16e-04) \\
\hline
\end{tabular}}
\end{table}

\begin{figure}[H]
\centering
\includegraphics[width=0.9\linewidth]{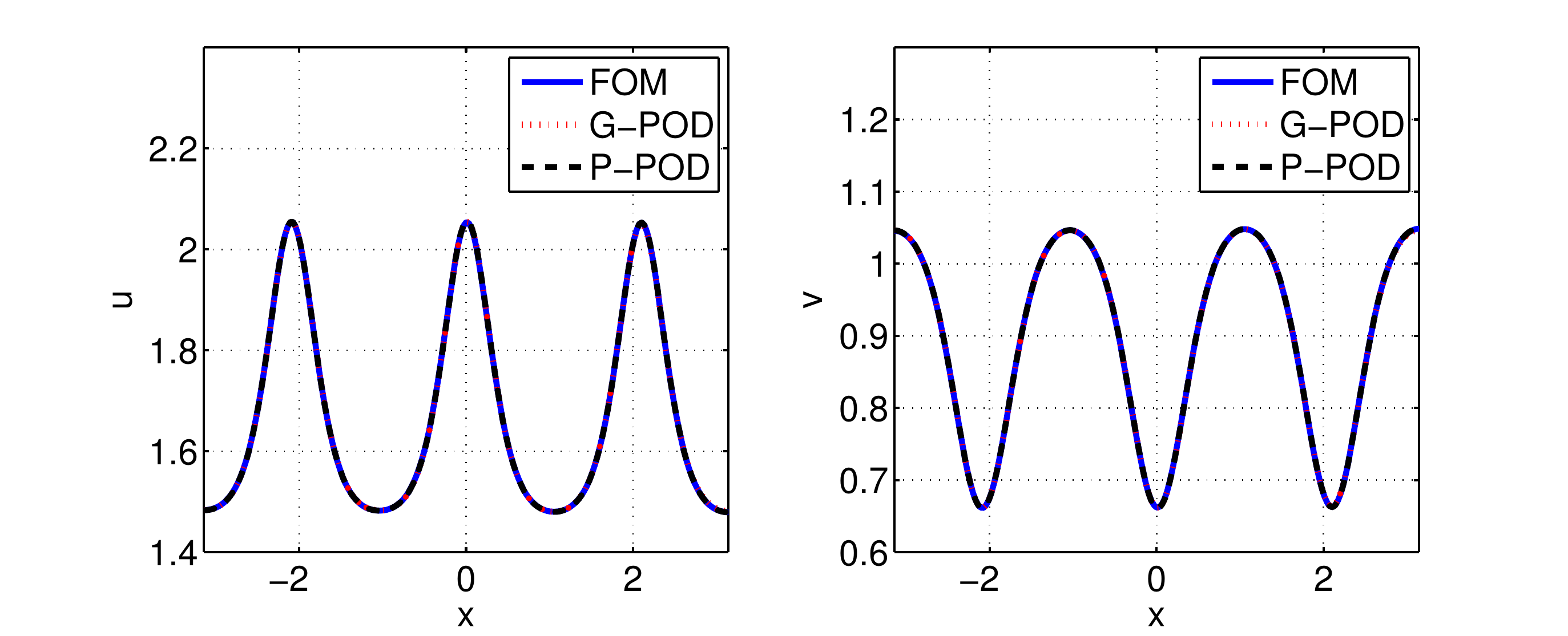}
\caption{FOM \& ROM  solutions for one-dimensional problem }
\label{1Dsol}
\end{figure}

A common way to test the quality of the ROMs depends on examining the ROMs beyond the regimes used for basis construction, which is generally called training interval. Moreover, the reliable predictions outside of the training interval show the amount of the dynamics captured by ROMs. In this section, to examine the quality of the ROM, we train the ROM until $ t=3 $ before it reaches the steady-state. In Figure~\ref{fig:den_pred1D}, we compare densities of the FOM solutions with the densities of the ROM solutions on testing and training interval for $ \text{tol}_\text{RIC} = 10^{-4} $. The densities in Figure~\ref{fig:den_pred1D} show that the ROM model is in good agreement with the FOM. Furthermore, we demonstrate the FOM and ROM solutions at time $ t=15 $ in Figure~\ref{fig:sol_pred1D}, which also indicates that the important dynamics of the FOM are captured.

\begin{figure}[H]
	\centering
    \includegraphics[width=0.45\linewidth]{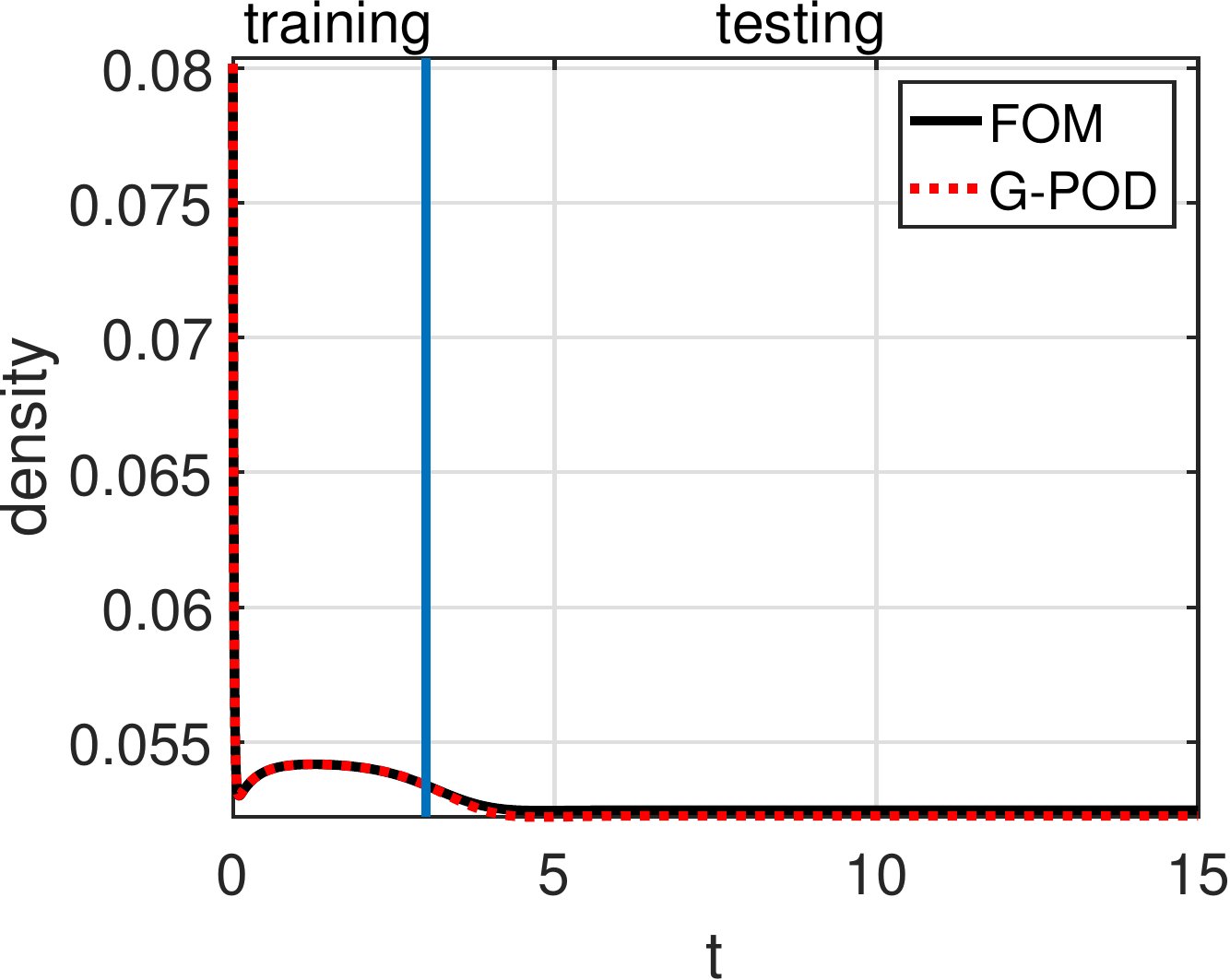}
	\includegraphics[width=0.45\linewidth]{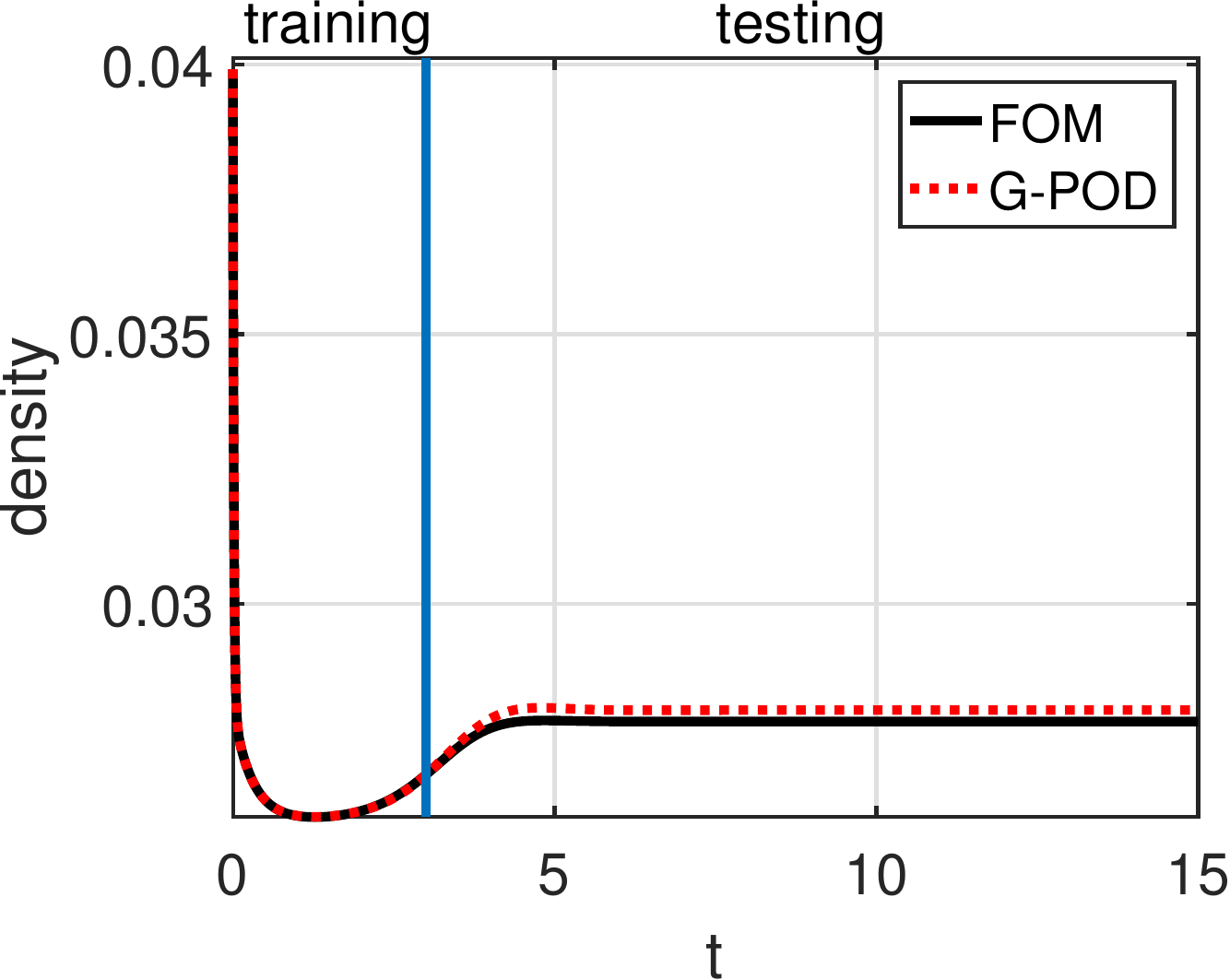}
	\caption{Predicted densities of $u$ (left) and $v$ (right) components for one-dimensional problem}
	\label{fig:den_pred1D}
\end{figure}

\begin{figure}[H]
	\centering
	\includegraphics[width=0.45\linewidth]{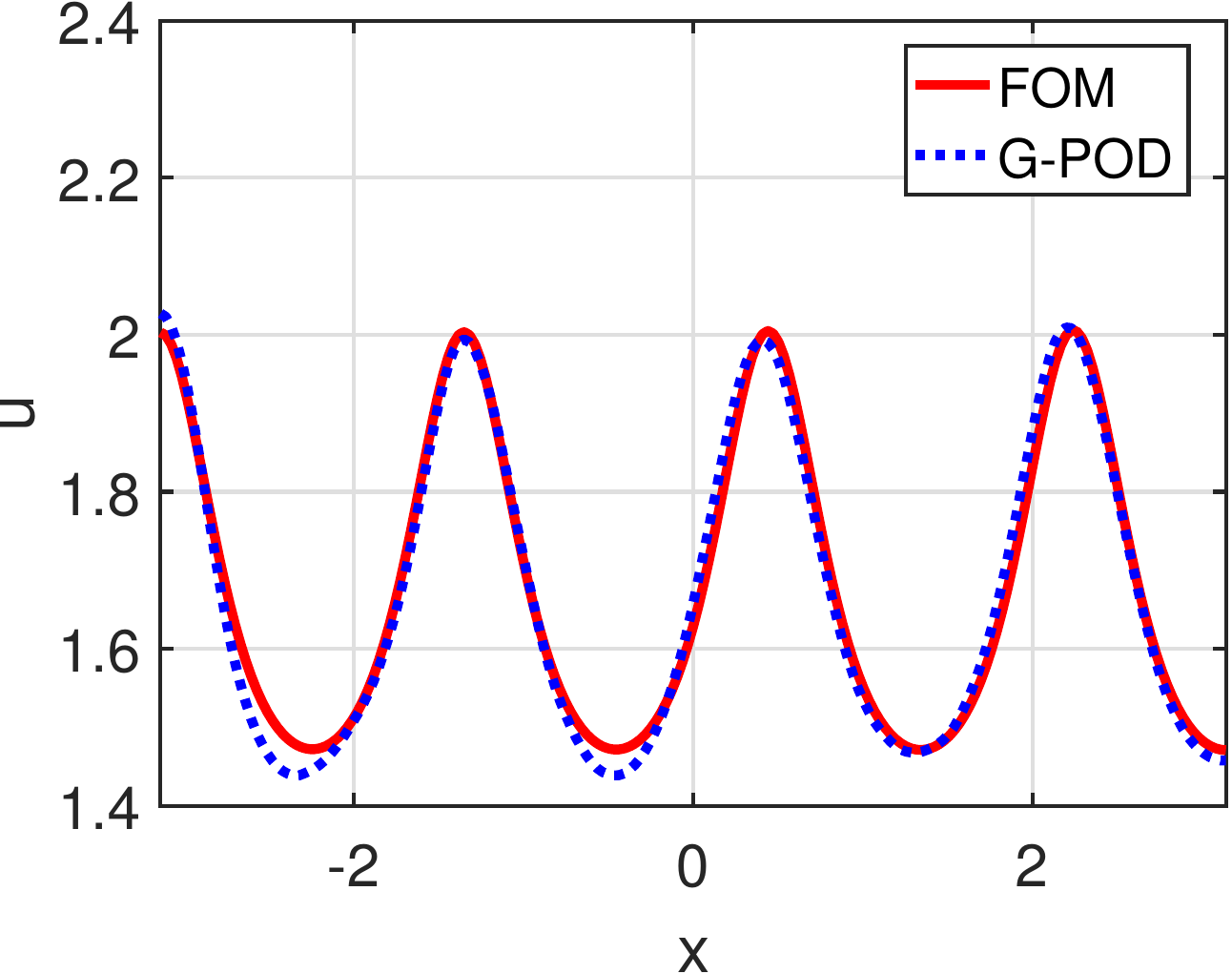}
	\includegraphics[width=0.45\linewidth]{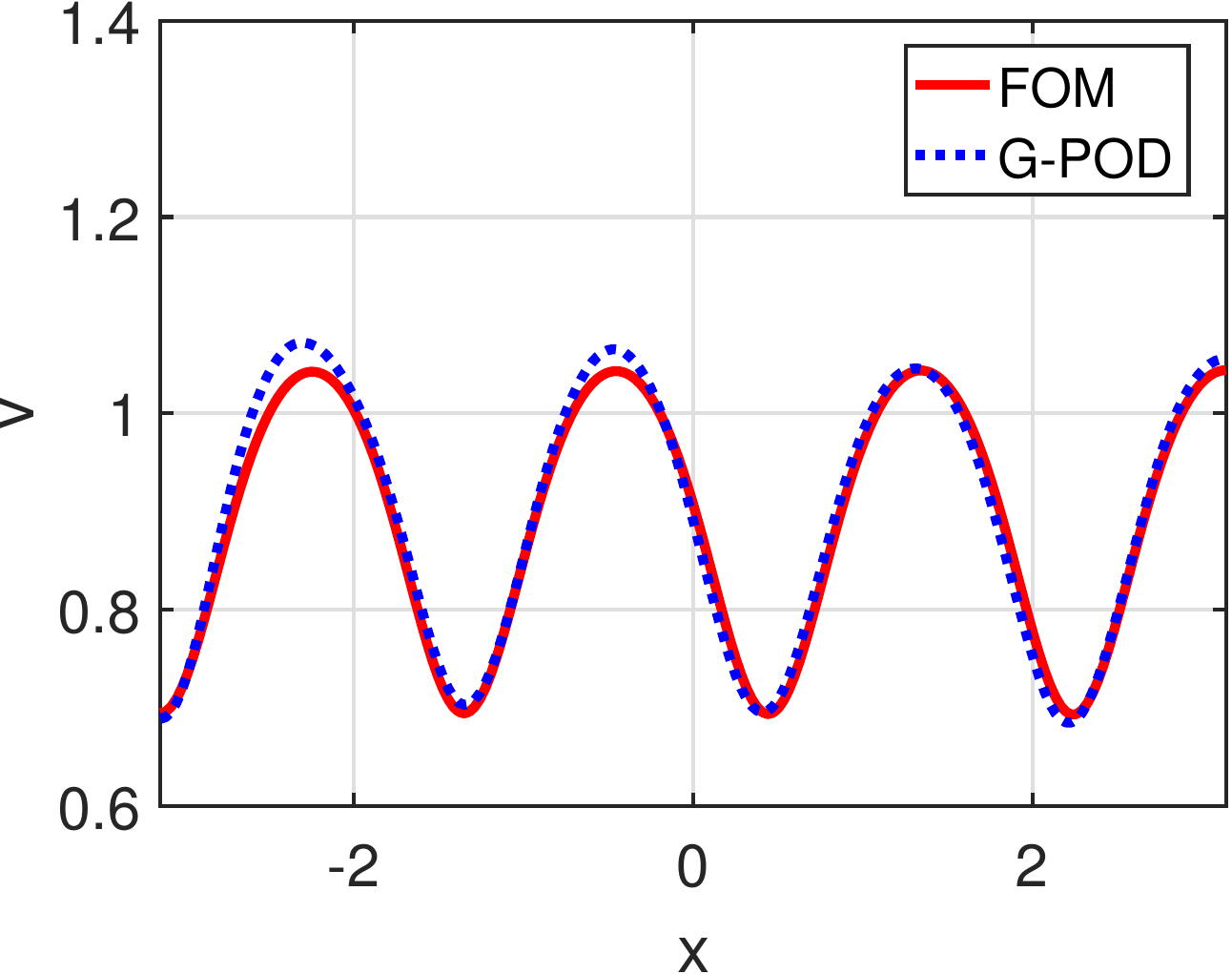}
	\caption{Predicted solutions of $u$ (left) and $v$ (right) components at $t=15$ for one-dimensional problem}
	\label{fig:sol_pred1D}
\end{figure}

\subsection{Two-dimensional SKT equation }
\label{exmp2}

Our second example is the two-dimensional SKT equation \eqref{skt} with
the parameters are taken from \cite{Gambino13}
\begin{align*}
a_1 &= 0.01,\; a_2 = 0.001, \; b_1 = 7.264, \; b_2 = 1.1, \; c_1 = 0.1,\; c_2 = 0.2, \\
\Gamma  &= 28.05, \;  r_1 =1.2,  \; r_2=1, \;\gamma_{11} = 0.5,\; \gamma_{12} = 0.4,\; \gamma_{21} = 0.38,\; \gamma_{22} = 0.4.
\end{align*}
Spatial domain is set to $\Omega =[0,\sqrt{2}\pi]\times [0,2\pi]$.  We take in both space directions the same number of partition $n_x=n_y=100$, and time step size is set to $\Delta t = 0.001$. The steady-states are reached at $t=2.938$.

The computational times for the matricized tensor $\widehat{H}= H(V_i\otimes V_j)\in\mathbb{R}^{N\times k_ik_j}$ versus the number of POD modes is plotted in Figure~\ref{cpu}
  using MULTIPROD, the MATLAB's loop over FOM dimension $N$ in \eqref{goyal} and  the Algorithm \eqref{algo1}. The computational efficiency of the MULTIPROD and the Algorithm \eqref{algo1} over  \eqref{goyal} with MATLAB's for loop can  be clearly seen. All the three approaches produce exactly the same ROM solutions.

\begin{figure}[H]
	\centering
	\includegraphics[width=0.5\linewidth]{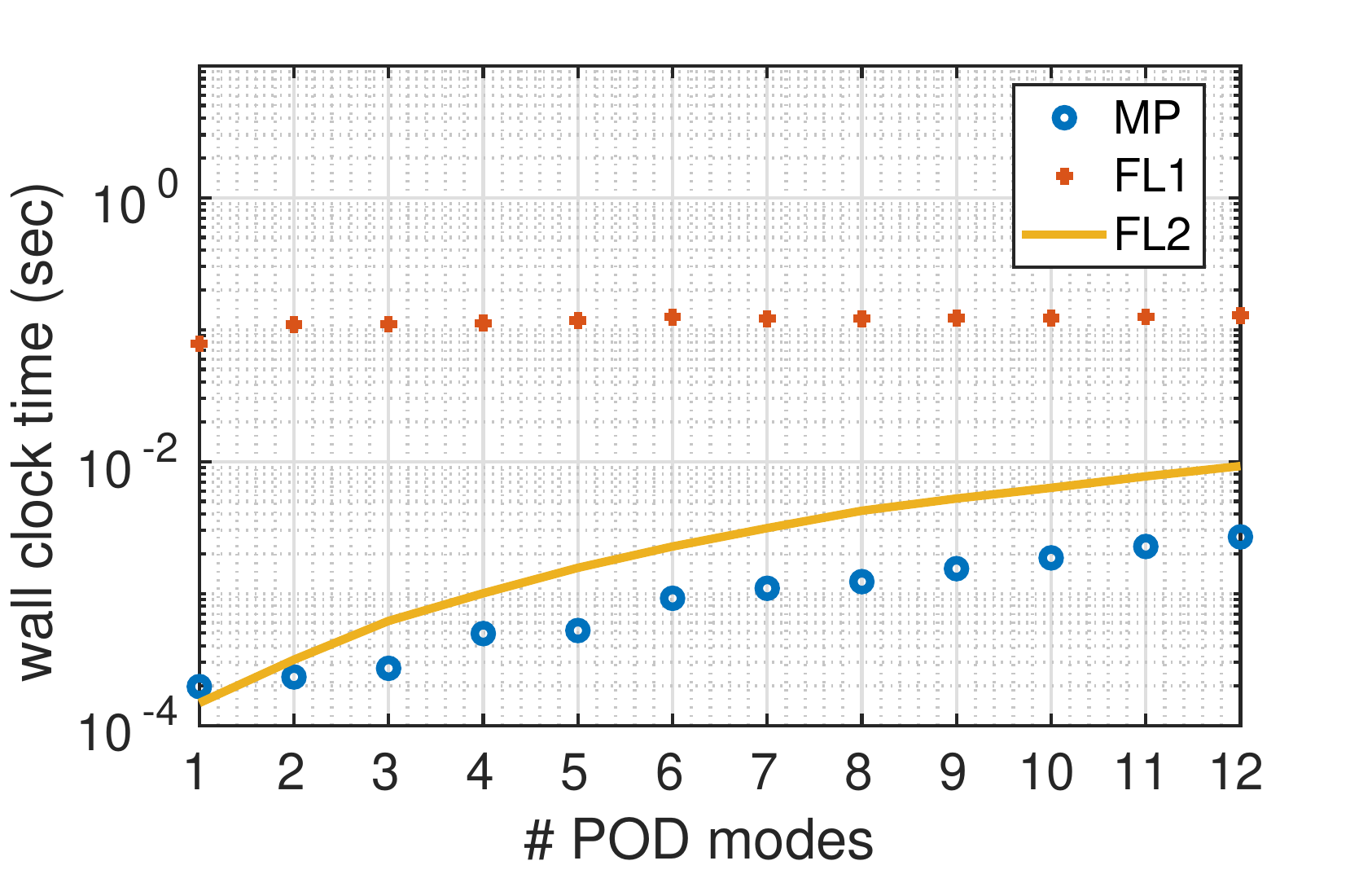}
	\caption{Wall clock times  for computing matricized tensor $\widehat{H}\in\mathbb{R}^{N\times k_ik_j}$ with MULTIPROD (MP), \eqref{goyal} (FL1) and the Algorithm \eqref{algo1} (FL2). }
	\label{cpu}
\end{figure}

In Figure~\ref{2Dden}, the densities start almost unchanged around $t=0.5$. By the PID tolerance $\text{tol}_{\text{PID}} = 10^{-7}$, the time interval is split into sub-intervals $I_1=[0, 0.484]$ and $I_2=[0.484, 2.938]$.
Decay of the singular values in Figure~\ref{2Dsvd} is similar to the one-dimensional SKT equation in the previous example. The FOM and ROM solutions in Figure~\ref{2Dfomsomsol} computed with the RIC tolerance $\text{tol}_{\text{RIC}}=10^{-6}$ agree well, and that the ones by the P-POD approach are almost the same as the FOM solutions. To show the time evolution  of the patterns by FOM and ROMs, a video component is also available for this example. To access the video components, simply click on one of the images in Figure~\ref{2Dfomsomsol} (online version only).
The errors by G-POD ad P-POD in Table~\ref{tab2D} show the same behavior as the errors in Table~\ref{tab1D} for the one-dimensional case.

\begin{figure}[H]
\centering
\includegraphics[width=0.45\linewidth]{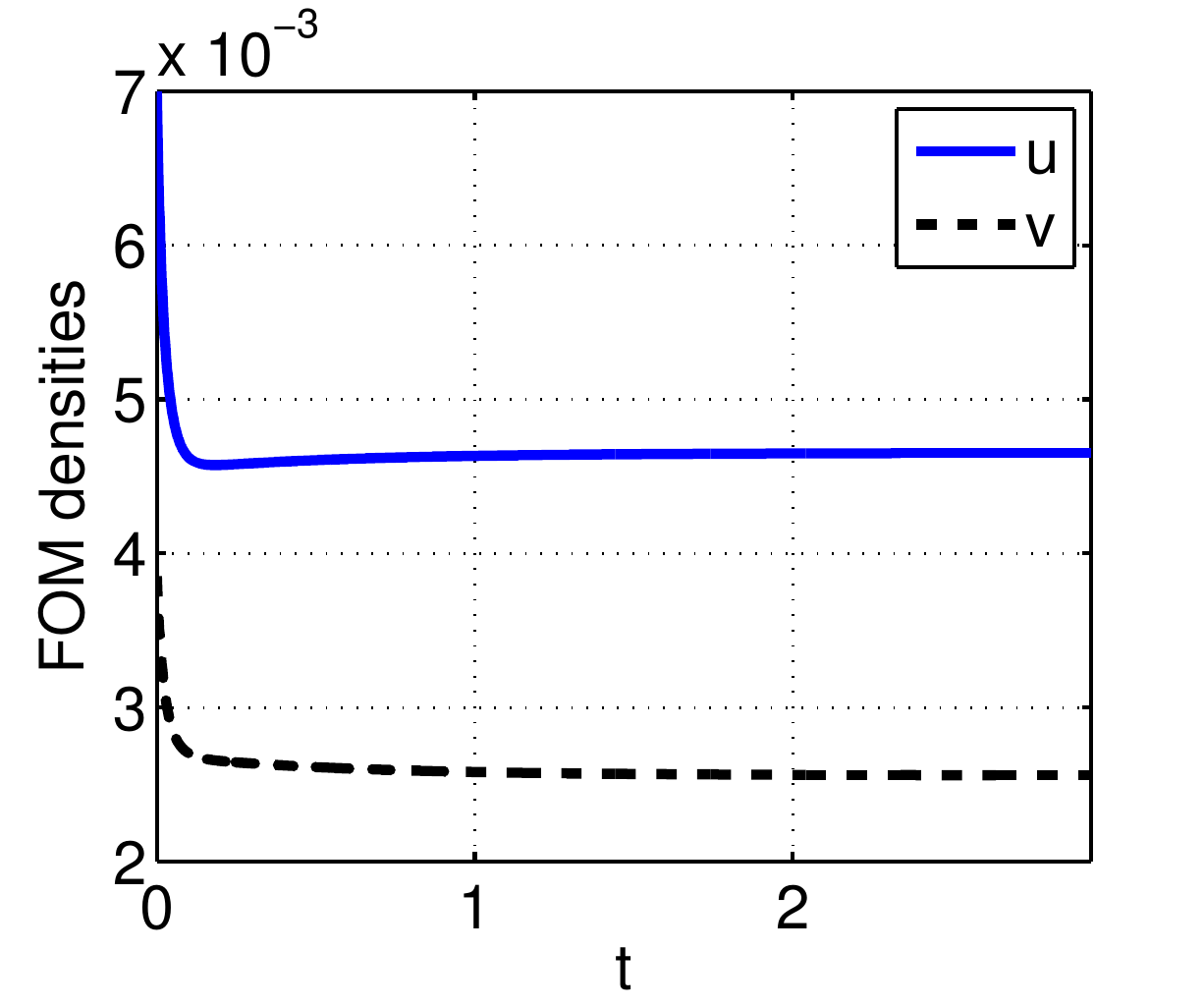}
\caption{FOM densities  for two-dimensional problem}
\label{2Dden}
\end{figure}

\begin{figure}[H]
\centering
\includegraphics[width=\linewidth]{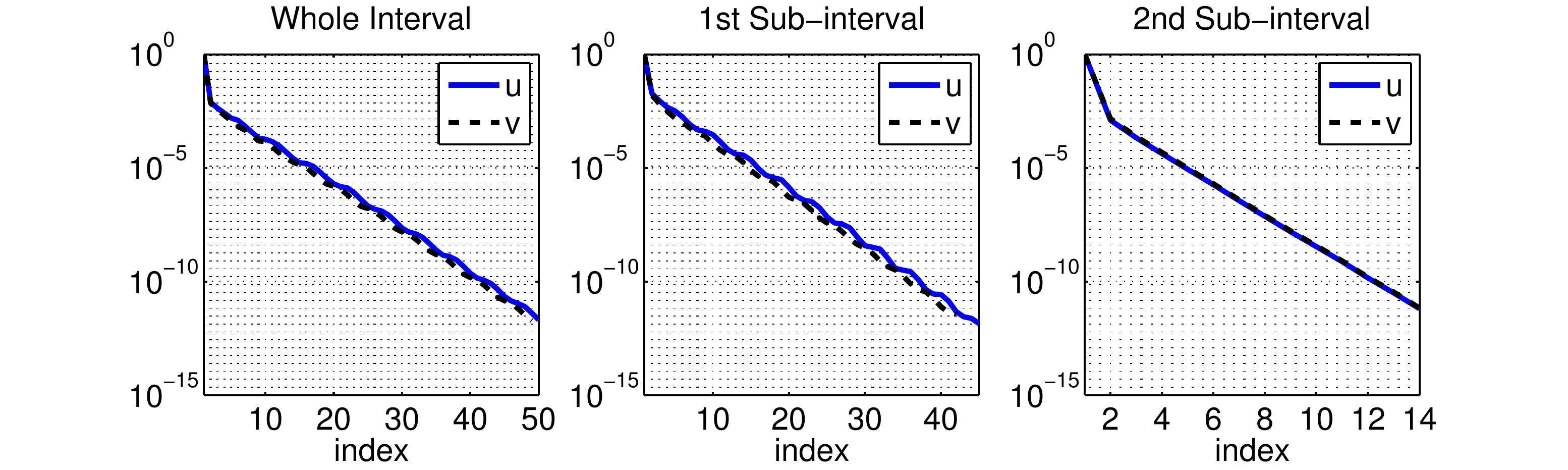}
\caption{Decay of normalized singular values  for two-dimensional problem}
\label{2Dsvd}
\end{figure}

\begin{figure}[H]
\centering
\includegraphics[width=\linewidth]{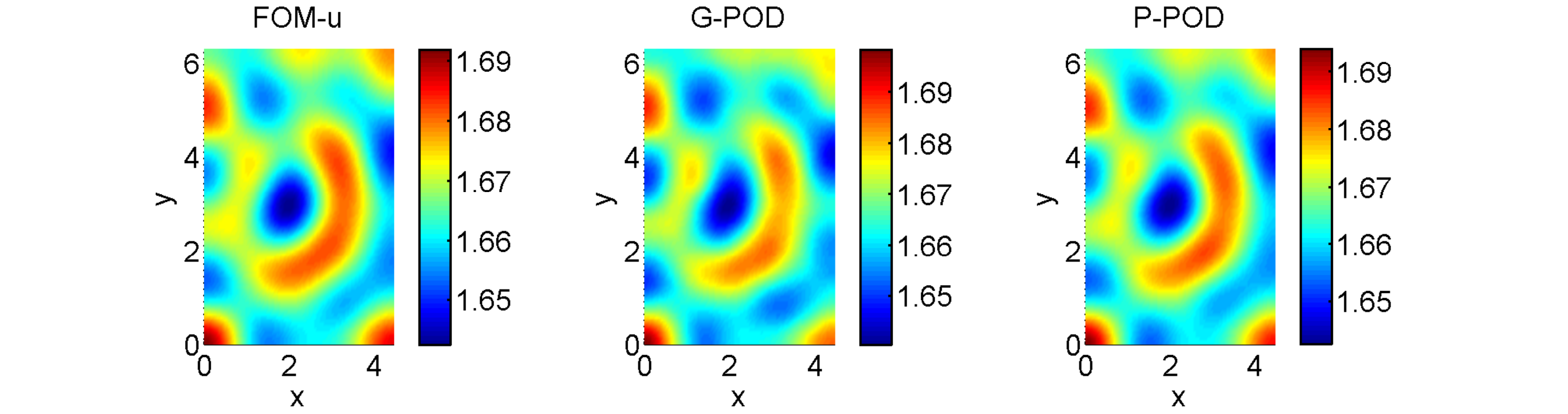}
\includegraphics[width=\linewidth]{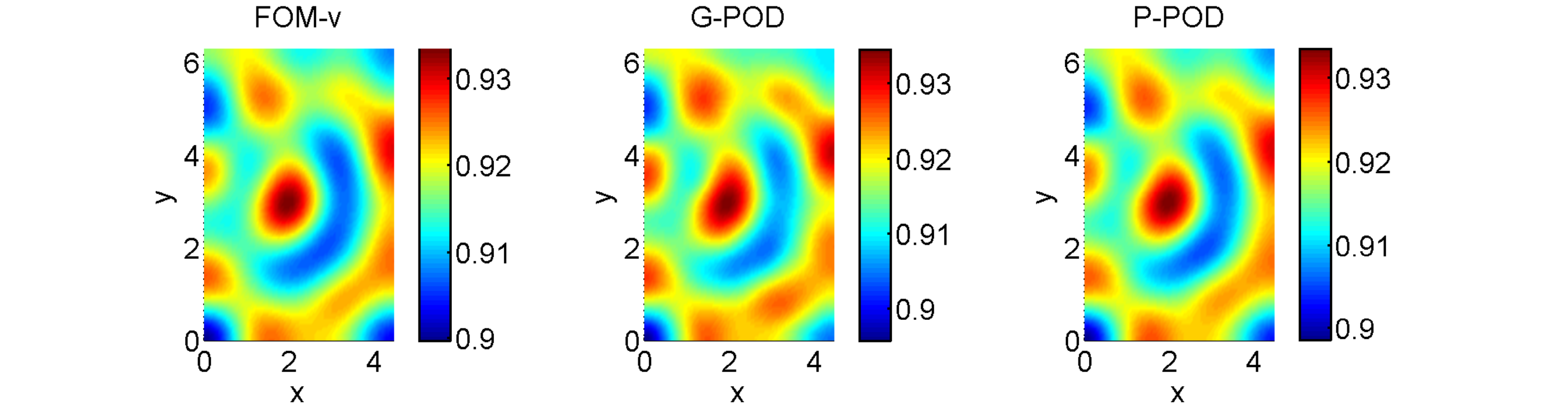}
\caption{FOM \& ROM  solutions of $u$ (top) and $v$ (bottom) components for two-dimensional problem }
\label{2Dfomsomsol}
\end{figure}

\begin{table}[H]
  \centering
	\caption{Number of POD modes and errors for two-dimensional problem \label{tab2D}}
\resizebox{\columnwidth}{!}{
\begin{tabular}{|c|ll|lll|}\hline
\multirow{2}{*}{$\text{tol}_\text{RIC}$} & \multicolumn{2}{c|}{G-POD $u(v)$} & \multicolumn{3}{c|}{P-POD $u(v)$} \\
 &  \#Modes  & Error   & \#Modes-$I_1$  & \#Modes-$I_2$  & Error  \\
\hline
$10^{-3}$ & 4(3) & 1.97e-03(2.18e-03) & 5(4) & 1(1) & 1.26e-03(1.40e-03) \\
$10^{-4}$ & 6(5) & 9.17e-04(1.02e-03) & 7(6) & 2(2) & 3.53e-04(3.89e-04) \\
$10^{-5}$ & 9(8) & 2.14e-04(2.30e-04) & 10(8) & 2(2) & 2.20e-04(2.46e-04) \\
$10^{-6}$ & 12(11) & 6.73e-05(7.29e-05) & 11(10) & 3(3) & 5.76e-05(6.42e-05) \\
\hline
\end{tabular}}
\end{table}

Similar to the previous 1D problem, we examine the quality of the ROM for the 2D problem by the prediction capabilities of the ROM. In this example, we train the ROM up to $ t=1 $ with $\text{tol}_\text{RIC}=10^{-5}$. Figure~\ref{fig:den_pred2D} shows that the ROM is reliably capable of predicting the densities up to the steady-state solution. Moreover, we show the solutions and the corresponding absolute errors at time $ t=3 $ in Figure~\ref{fig:sol_pred2D}, which again indicates that the ROM captures all the important dynamics of the FOMs.

\begin{figure}[H]
	\centering
	\includegraphics[width=0.45\linewidth]{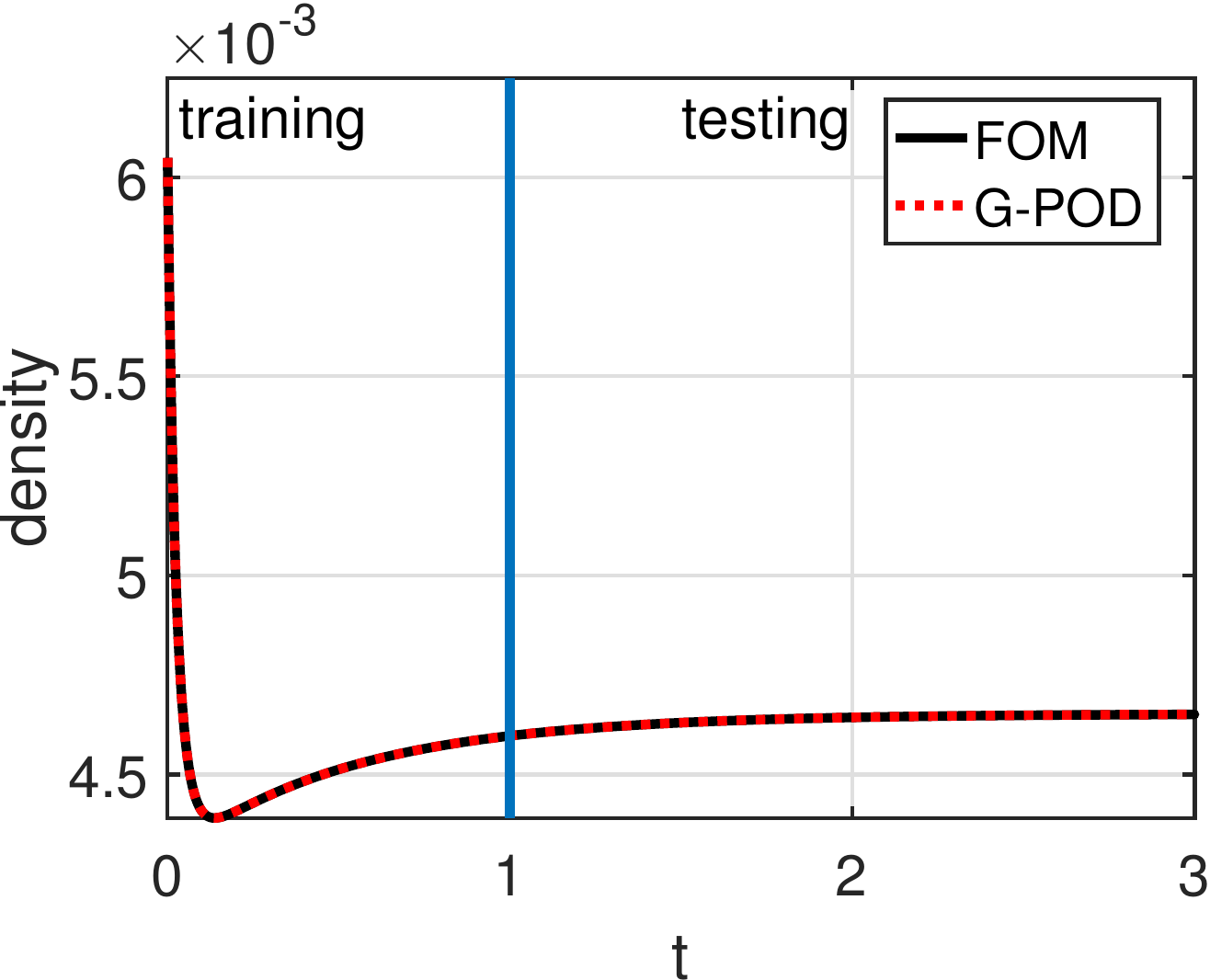}
	\includegraphics[width=0.45\linewidth]{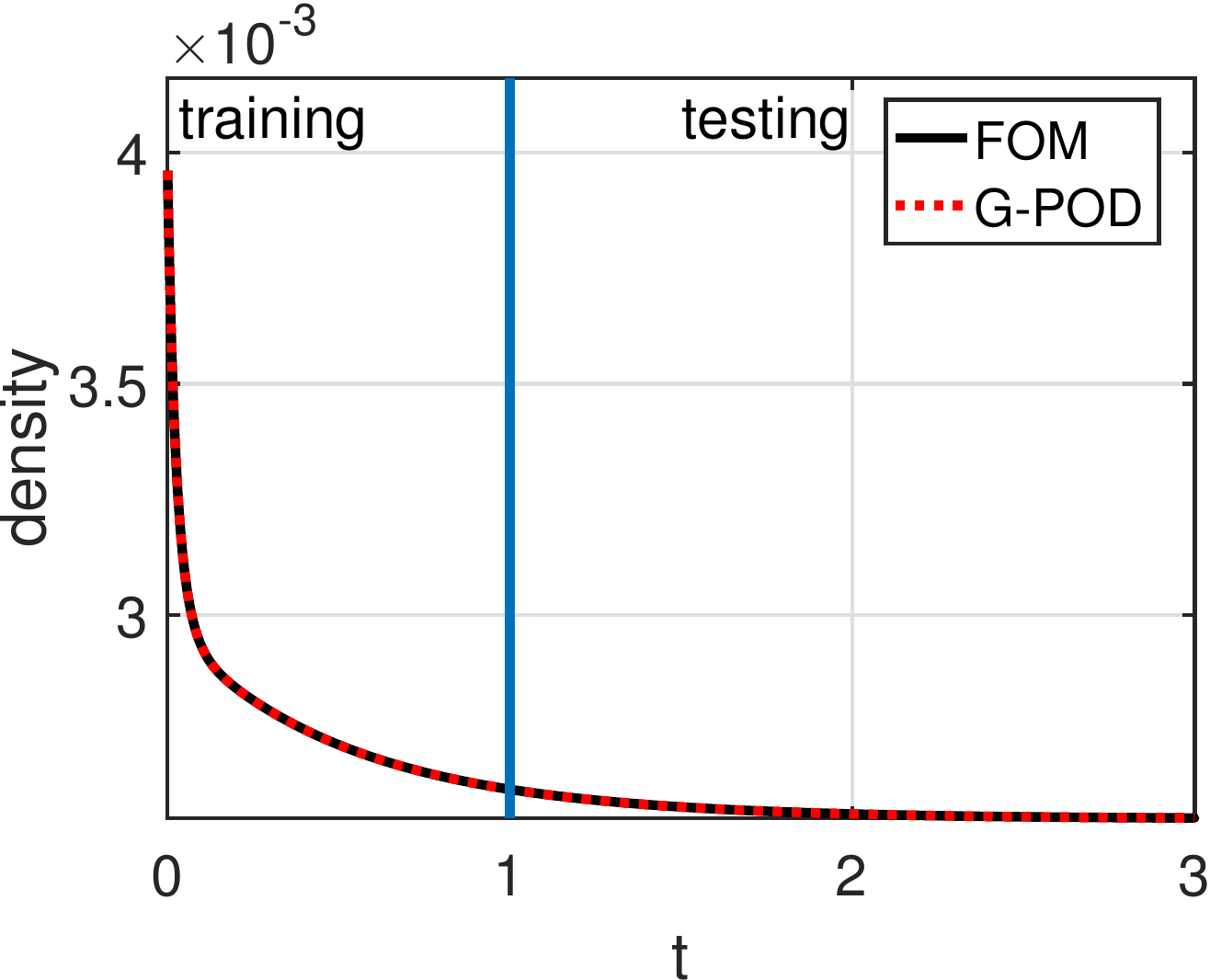}
	\caption{Predicted densities of $u$ (left) and $v$ (right) components for two-dimensional problem}
	\label{fig:den_pred2D}
\end{figure}

\begin{figure}[H]
	\centering
	\includegraphics[width=0.3\linewidth]{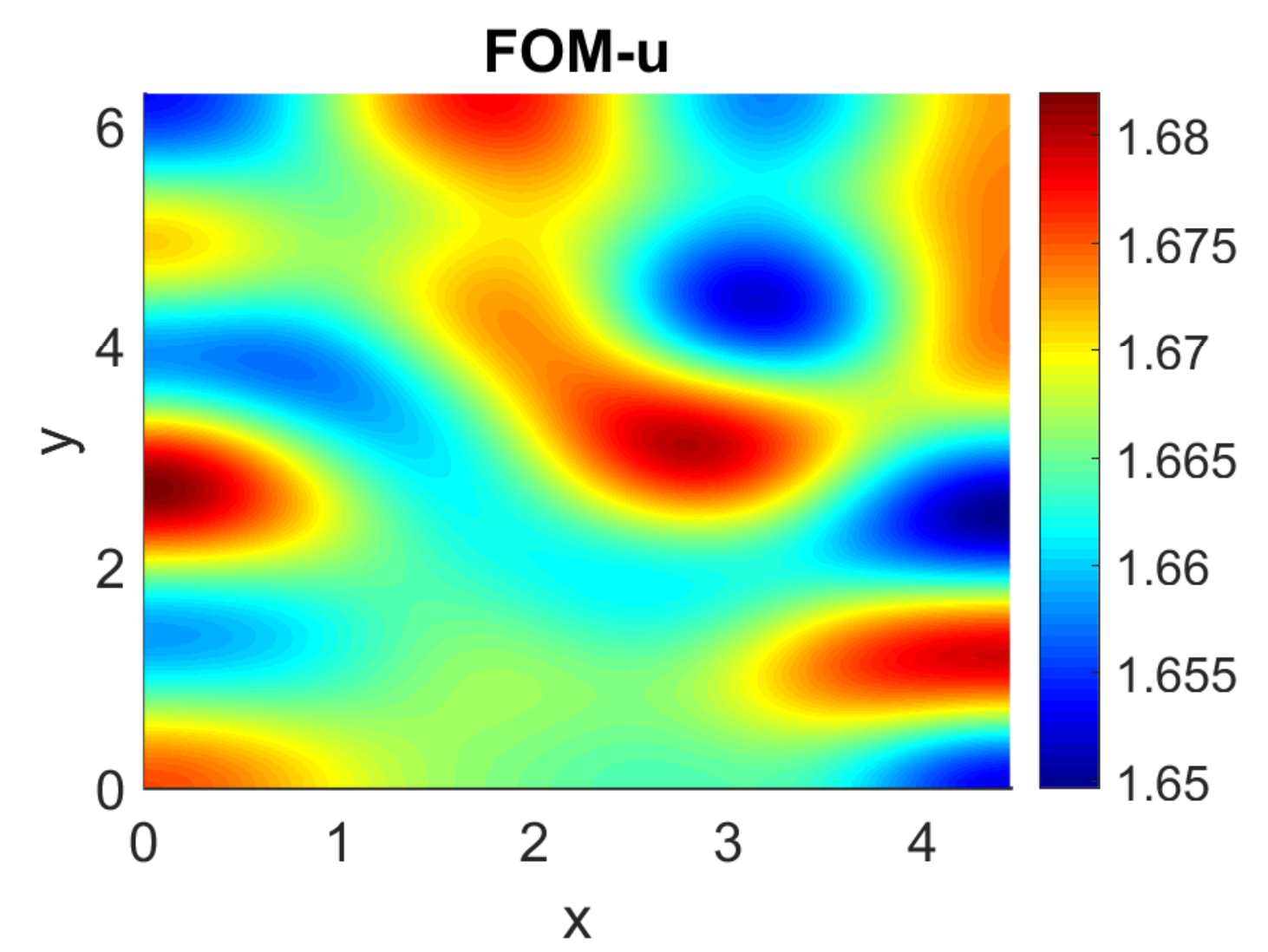}
    \includegraphics[width=0.3\linewidth]{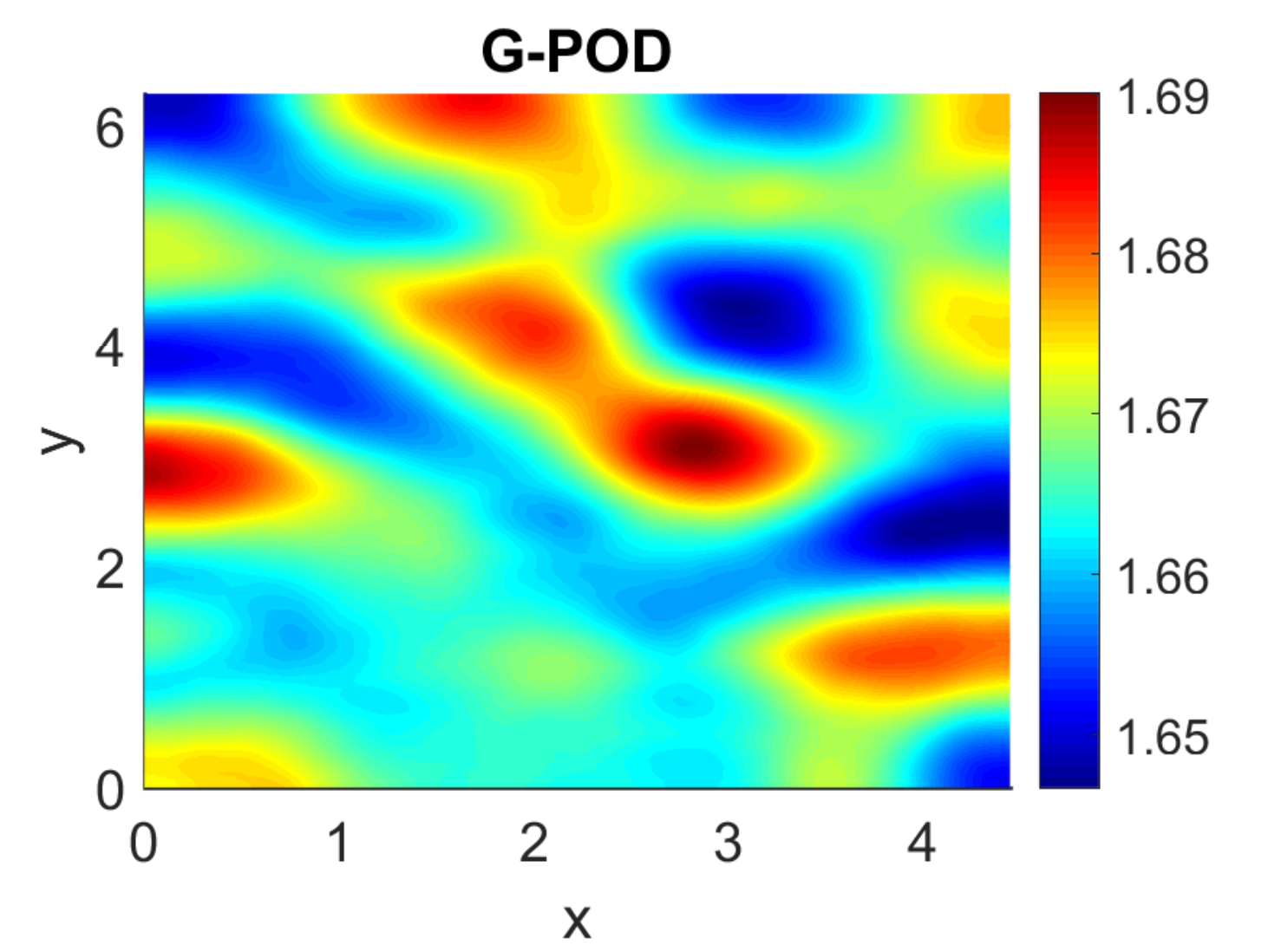}
    \includegraphics[width=0.3\linewidth]{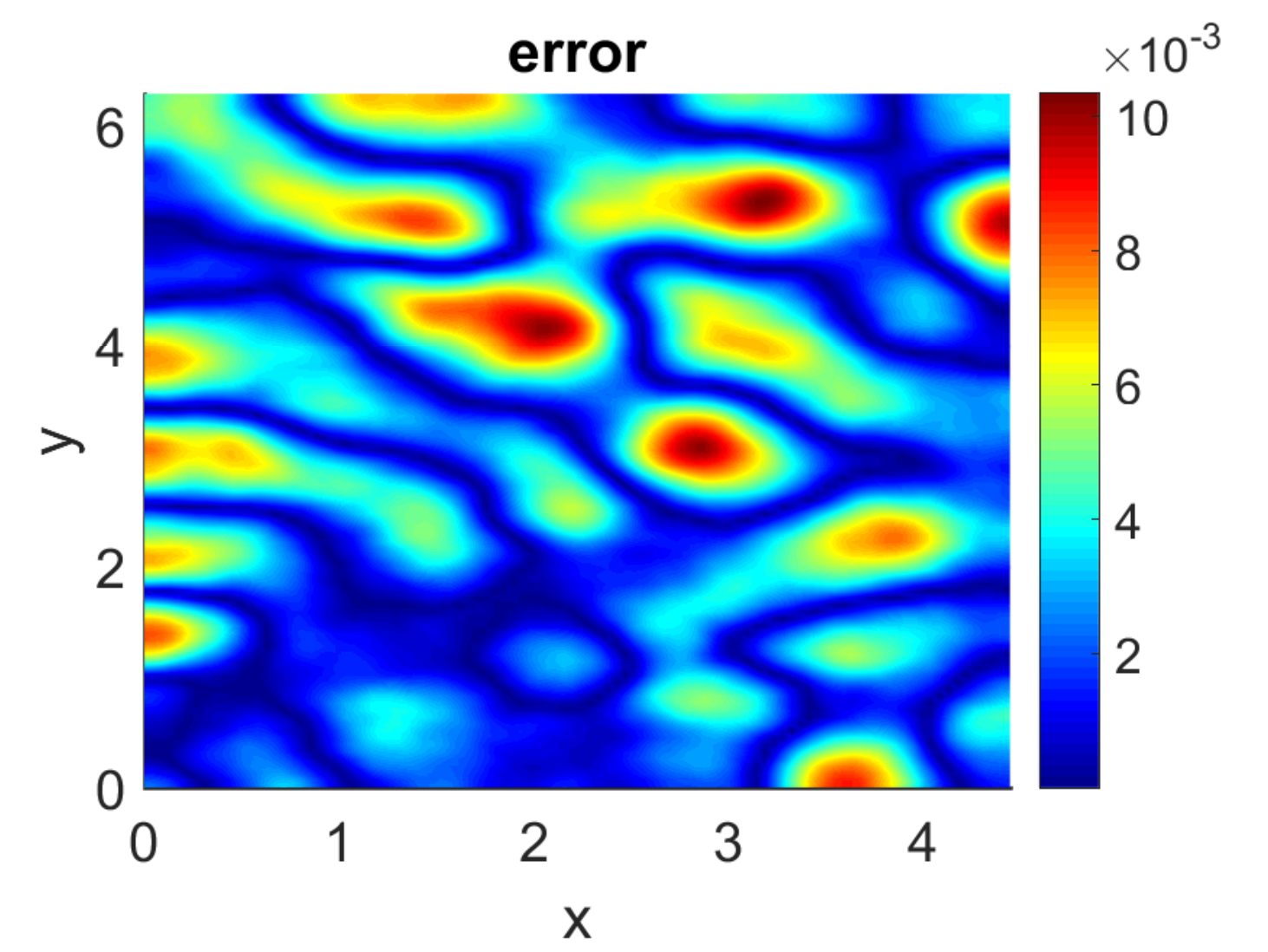}
    \includegraphics[width=0.3\linewidth]{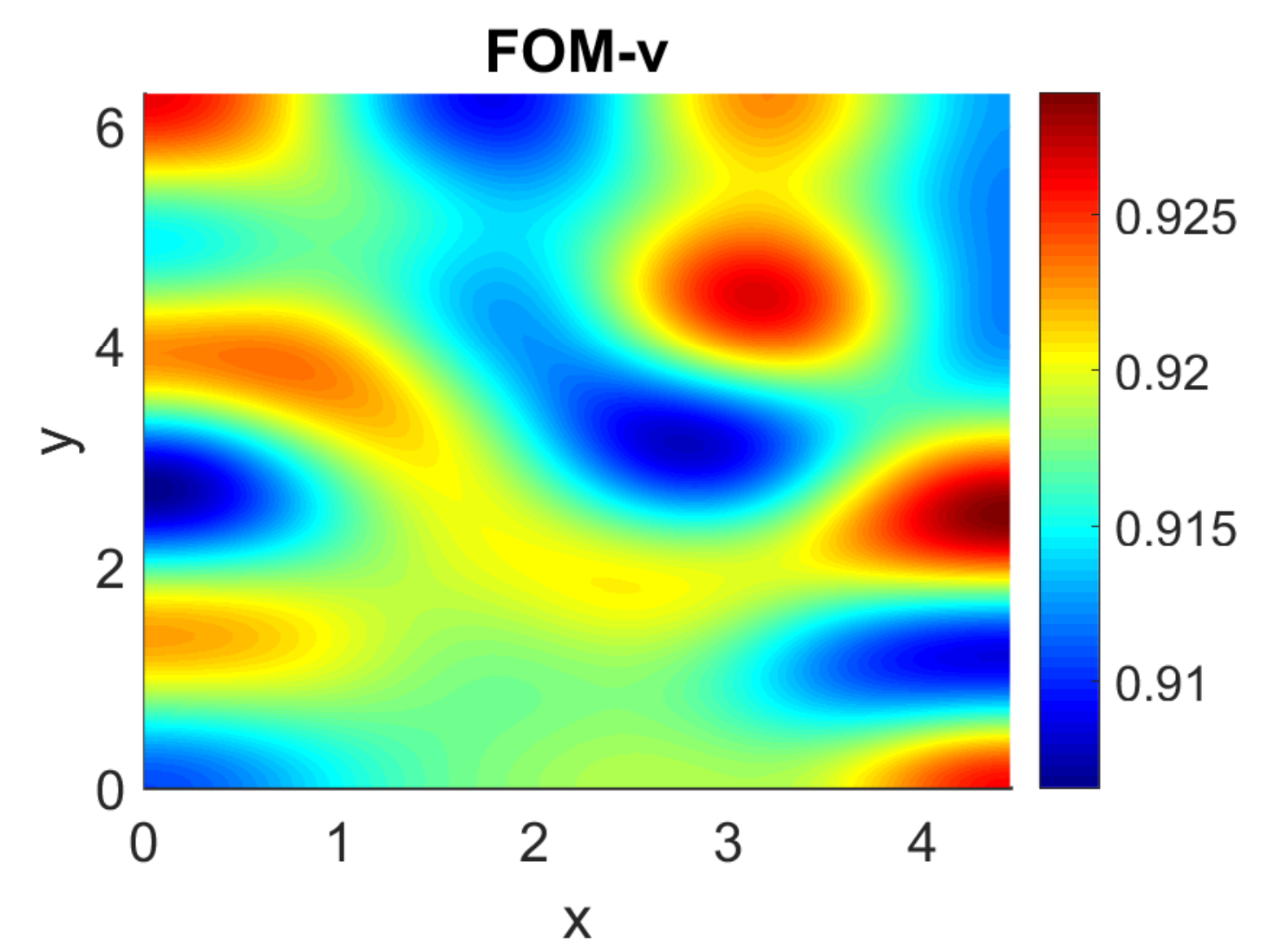}
    \includegraphics[width=0.3\linewidth]{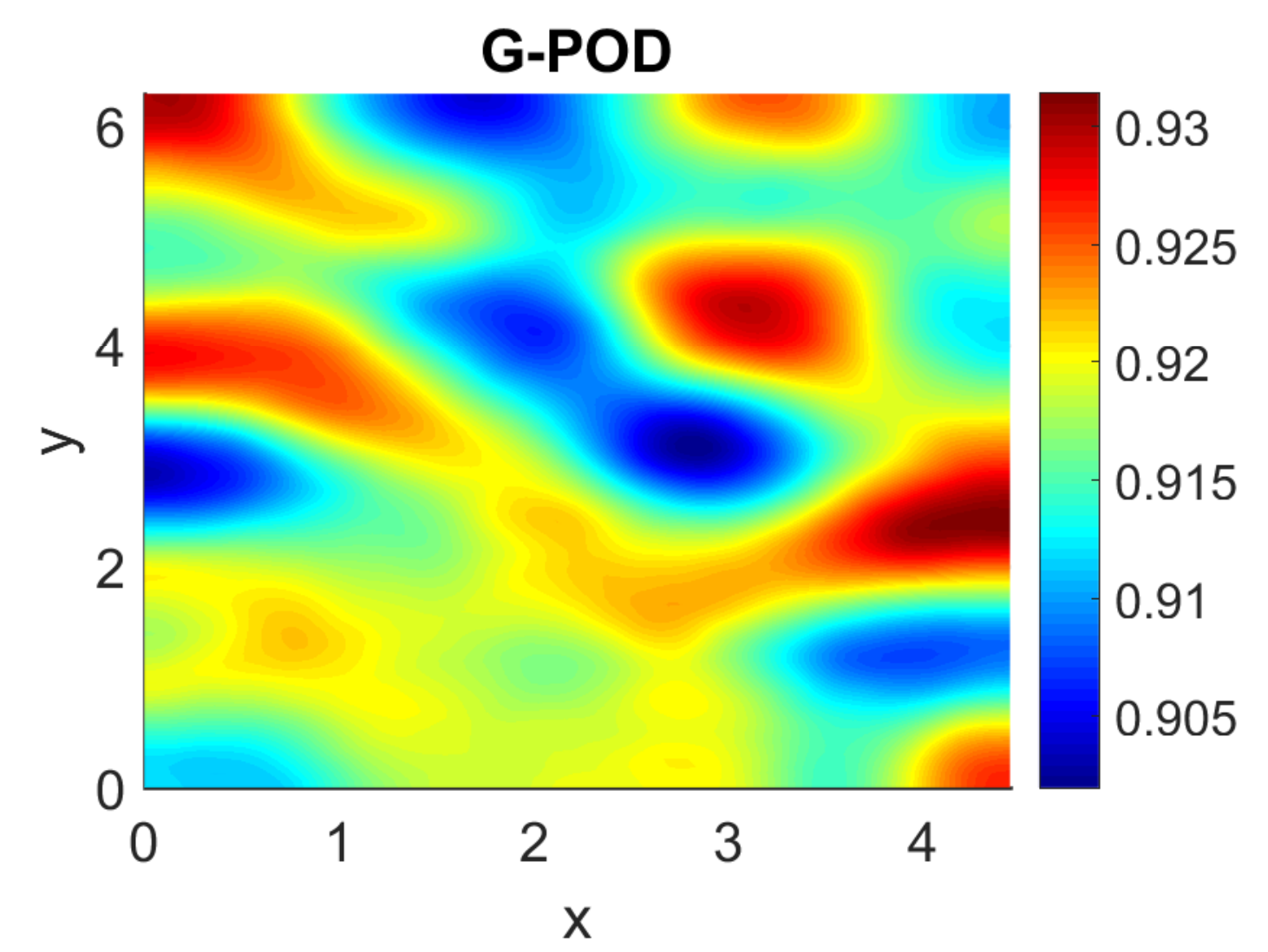}
    \includegraphics[width=0.3\linewidth]{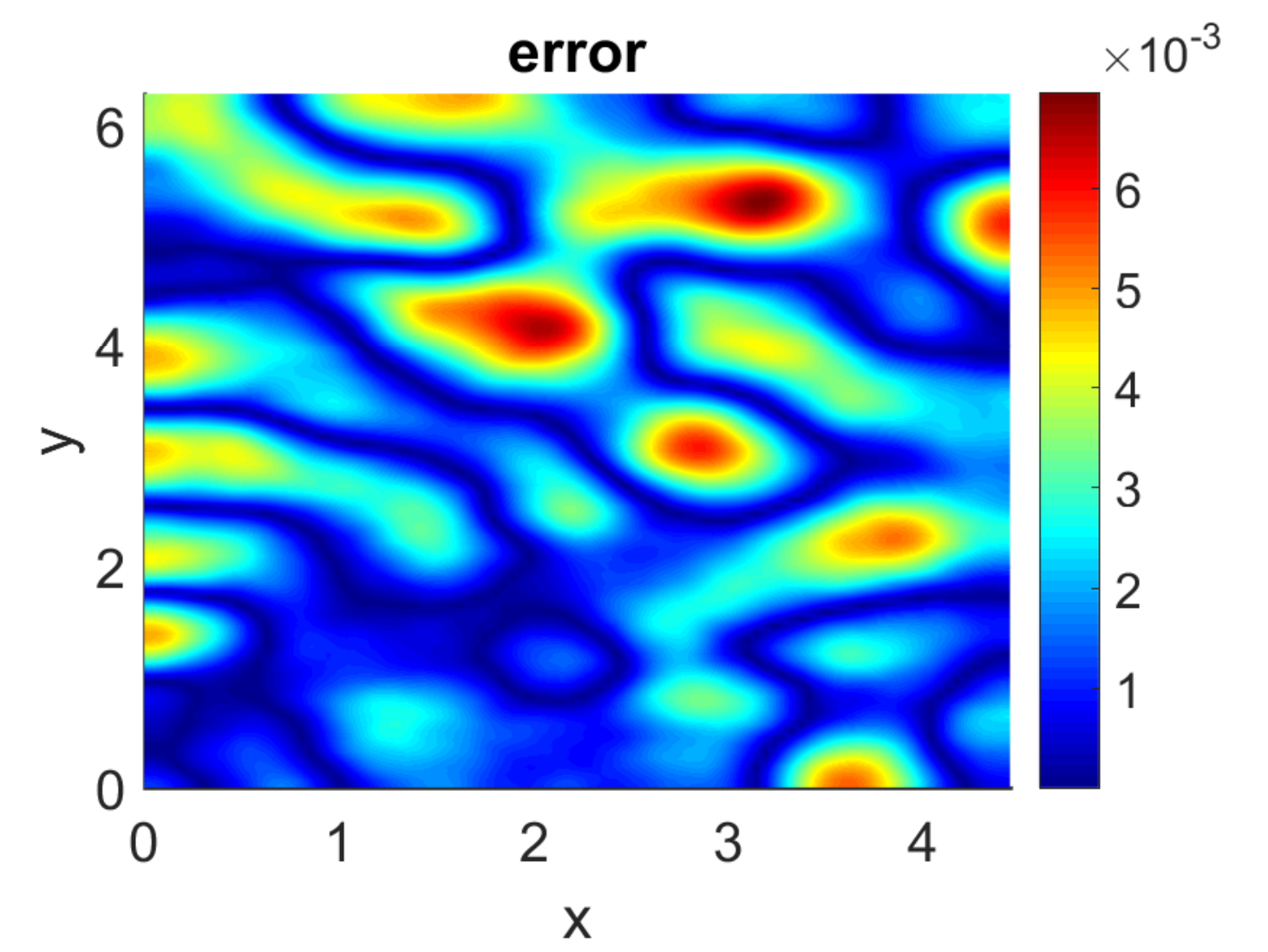}
	\caption{Predicted solutions of $u$ (top) and $v$ (bottom) components at $t=3$ for two-dimensional problem}
	\label{fig:sol_pred2D}
\end{figure}

\subsection{Entropy  preservation}

The entropy $\mathcal E$ in \eqref{ent} of the SKT equation is defined with the Lotka-Volterra kinetics  terms $f_i(u) =0$, $i=1,2$. The entropy is computed with the same diffusion coefficients  for one- and two-dimensional SKT equation \eqref{skt}. Initial conditions are given by \cite{Sun19}
\begin{align*}
u_1(x) & =  e^{\frac{1}{2}\sin x}, \quad u_2(x) = e^{\frac{1}{2}\cos  2x }, \\
u_1(x,y) & = \frac{1}{2} (\sin (\pi(x+y)))+ 1, \quad u_2(x,y) = \frac{1}{2} (\cos (\pi(x-y)))+1,
\end{align*}
for one- and two-dimensional problems, respectively. Since the SKT equation is solved without the reaction terms, the transient phase is absent. Therefore, the ROMs are computed only by the G-POD approach.
Figures~\ref{1Dent}-\ref{2Dent} show that the FOM entropy for both one- and two-dimensional problems decay with the time, and the same structure is well preserved by the entropy obtained by ROM.

\begin{figure}[H]
\centering
\includegraphics[width=0.9\linewidth]{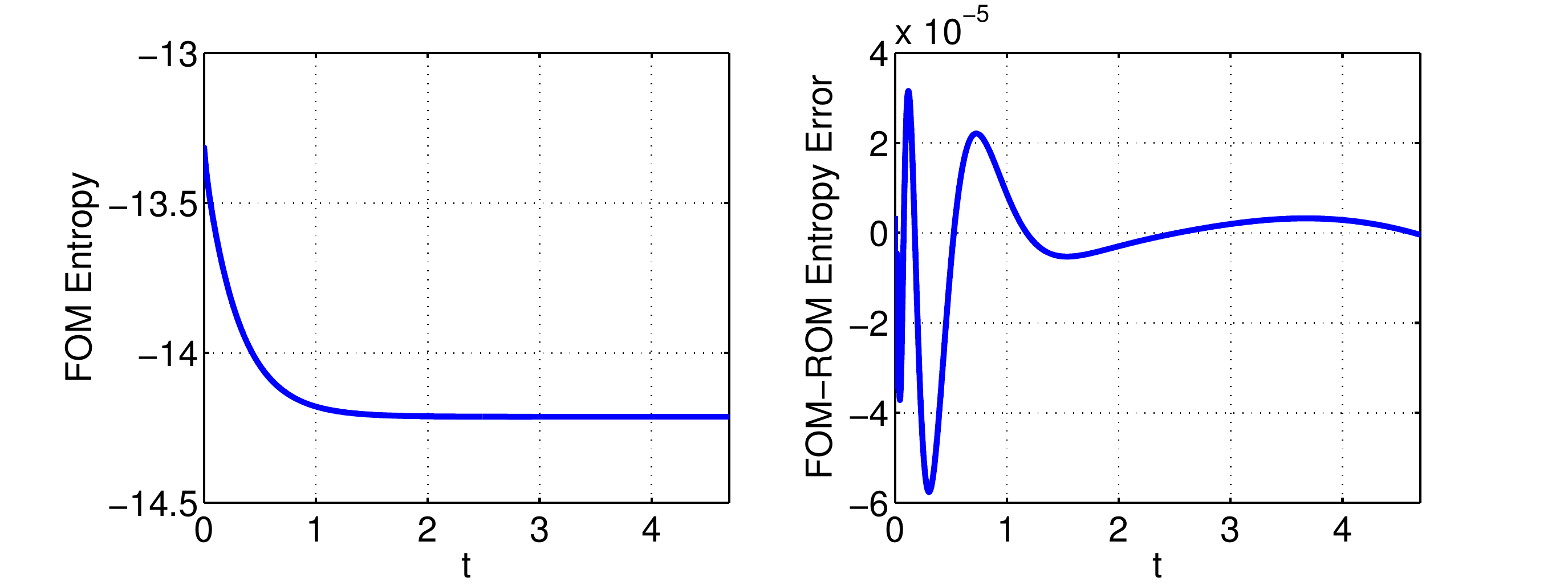}
\caption{Entropy decay for one-dimensional problem }
\label{1Dent}
\end{figure}

\begin{figure}[H]
\centering
\includegraphics[width=0.9\linewidth]{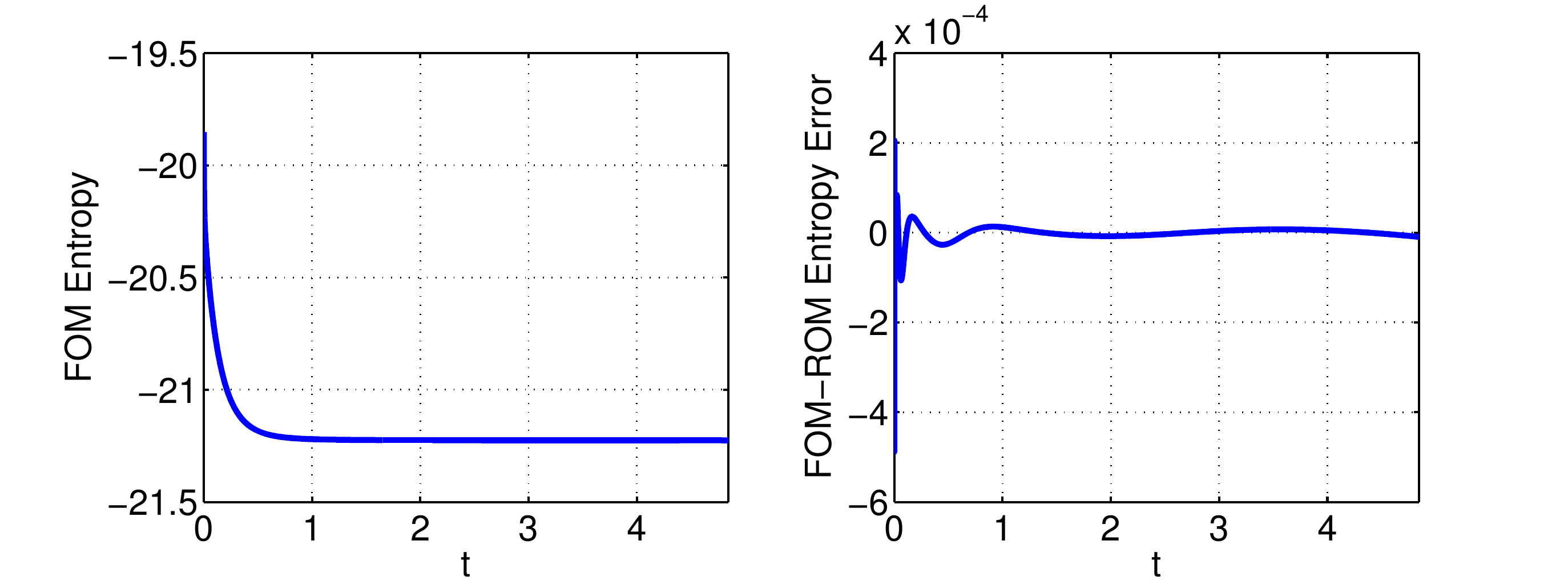}
\caption{Entropy decay for two-dimensional problem }
\label{2Dent}
\end{figure}

\section{Conclusions}

Exploiting the different behavior of transient and steady-state solutions of the SKT equation, reduced solutions are obtained in a computationally efficient way. The quadratic nonlinear terms of SKT equation are reflected in the semi-discrete linear-quadratic ODE system using finite-differences, which enables separation of the offline-online computation. The ROM solutions depend affinely on the parameters in both of the linear and quadratic parts. This allows the prediction of patterns for different parameter values without interpolation. We plan to investigate the bifurcation behavior of the SKT equation \cite{Kuehn20,Soresina19} using ROM techniques.

\noindent{\bf Acknowledgemets\/}
The authors thank for the constructive comments of the referees, which helped much to improve the paper.


\end{document}